%% file: main.tex
\documentclass{article}

\usepackage[english]{babel}

\usepackage[fleqn]{amsmath}
\usepackage{xcolor}
\usepackage{etoolbox} %
\definecolor{mypink}{RGB}{0, 0, 0}

\newcommand*\linenomathpatch[1]{%
  \cspreto{#1}{\linenomath}%
  \cspreto{#1*}{\linenomath}%
  \csappto{end#1}{\endlinenomath}%
  \csappto{end#1*}{\endlinenomath}%
}

\linenomathpatch{equation}
\linenomathpatch{gather}
\linenomathpatch{multline}
\linenomathpatch{align}
\linenomathpatch{alignat}
\linenomathpatch{flalign}

\usepackage{amsfonts}
\usepackage{amssymb}
\usepackage{amsthm}
\usepackage{centernot}
\usepackage{graphicx}
\usepackage{subcaption}
\usepackage{tabularx}
\usepackage{multirow}
\usepackage{algorithmic}
\usepackage[ruled,vlined]{algorithm2e}
\newtheorem{thm}{Theorem}[section]

\newtheorem{cor}[thm]{Corollary}

\newtheorem{defn}[thm]{Definition}

\usepackage[colorlinks=true, allcolors=blue]{hyperref}
\DeclareMathOperator{\tr}{tr}

\usepackage{fullpage}
\usepackage{lineno}

\usepackage{url}
\newcommand{\bm}{\mathbf}

\usepackage{xcolor}
\usepackage{authblk}

\title{Probing robustness of nonlinear filter stability numerically using Sinkhorn divergence}
\author[1]{Pinak Mandal\thanks{Corresponding author: \texttt{pinak.mandal@icts.res.in}}}
\author[1]{Shashank Kumar Roy}%
\author[1, 2]{Amit Apte}%

\affil[1]{{\small International Centre for Theoretical Sciences TIFR, Bangalore 560089 India}}
\affil[2]{{\small Indian Institute of Science Education and Research, Pune 411008 India}}

\begin{document}

\maketitle

\begin{abstract}

Using the recently developed Sinkhorn algorithm for approximating the Wasserstein distance between probability distributions represented by Monte Carlo samples, we demonstrate exponential filter stability of two commonly used nonlinear filtering algorithms, namely, the particle filter and the ensemble Kalman filter, for deterministic dynamical systems. We also establish numerically a relation between filter stability and filter convergence by showing that the Wasserstein distance between filters with two different initial conditions is proportional to the bias or the RMSE of the filter.

\end{abstract}

\section{Introduction}\label{sec-intro}

\input{01-intro}

\section{Problem statement}\label{sec-problem}
\input{02-problem}

\section{Methodology}\label{sec-method}
\input{03-method}

\section{Results}\label{sec-results}
\input{045-filter-results}

\section{Conclusions}\label{sec-conclude}

The main focus of this paper is the numerical study of nonlinear filter stability. For this purpose, we use the recently developed Sinkhorn algorithm to calculate an approximation of the Wasserstein distance between Monte Carlo samples from probability distributions. This allows us to directly assess the stability by computing the expected value, averaged over multiple observation realizations, of the distance between filtering distributions as function of time. We also studied extensively the dependence of stability properties on two main parameters, the time between observations and the observational error covariance. 

Our results provide strong numerical indication of exponential stability of PF and EnKF for deterministic chaotic dynamical systems. For a fixed observational covariance, as the gap between the observations increases, the exponential rate of decay of the distance $D_\varepsilon$ between two filters decreases for the particle filter whereas it remains approximately constant and close to twice the Lyapunov exponent for the ensemble Kalman filter. {\color{mypink} Further exploration of relationships between the chaotic properties of the system and the exponential stability is an interesting direction for future work.} Mathematical proofs of stability of these numerical algorithms is another interesting open area of further research.

Further, with increasing observational gap, the filter uncertainty and bias decrease for PF but increases for the EnKF. In general, the EnKF has significantly smaller uncertainty and bias as compared to the PF. More extensive numerical explorations and associated theoretical studies to understand the differences between these two important numerical filters - the particle filter and the ensemble Kalman filter - is an important direction for further investigations.

The other main focus was to establish a relation between the filter stability and filter convergence. We show that for a wide variety of parameters, the distance between two filters started with different initial conditions is directly proportional to the $l_2$ error or the RMSE between the filter mean and the true underlying trajectory. The significance of this result is that the techniques developed in this paper provide tools to check filter stability, with synthetic or real observations for a given dynamical and observational model, even before starting to use a filter with actual observations, while the results connecting stability with convergence of the $l_2$ error can then be used to give indication of whether the filter is also accurate. This is particularly useful in real-world scenarios where one does not have access to the true trajectory of the dynamical system under consideration and hence can not compute the RMSE. In these cases one can still perform assimilation with two different initial conditions and compute the distance between the resulting filtering distributions and how this distance changes in time can indicate how the accuracy of the filter changes with time since these two quantities are strongly correlated. 

\section*{Acknowledgements}
This work was supported by the Department of Atomic Energy, Government of India, under project no. RTI4001. The authors would like to thank the anonymous referees for the comments which helped improve the manuscript.

\section{Appendix}\label{sec-append}
\input{06-appendix}

\bibliographystyle{alpha}
\bibliography{ref}
\end{document}

%% file: 01-intro.tex
The problem of predicting the state of a complex dynamical system is
ubiquitous in many scientific and engineering fields. In the context
of earth sciences, weather prediction and reanalysis of past climate
are major examples of such state estimation problems, which have two
main ingredients: (i) a dynamical model, usually deterministic, of the
system and (ii) partial, usually sparse, and needless to say, noisy,
observations. The process of combining these observations with the
model to get an ``optimal'' state estimate is commonly called data
assimilation - a term introduced in earth sciences~\cite{Carrassi-etal2008, fletcher2017data, Kalnay03}.

The mathematical formulation of this problem in a Bayesian framework
encapsulates the information from the model in terms of a prior
distribution, and the observational likelihood is used to obtain a
posterior distribution for the model state~\cite{ApteH07, law2015data}. For dynamical systems,
this is precisely the problem of nonlinear filtering, where the
posterior distribution itself changes with time and is conditioned on
observations up to that time. This posterior is called the {\it
filter} or the {\it filtering distribution}~\cite{sarkka2013bayesian, vanleeuwen2015book, asch2016data}.

A crucial characteristic of the atmospheric and oceanic dynamics is
their chaotic nature, manifested in sensitive to initial conditions. Thus a natural question is
whether the filter is also sensitive to the choice to the initial
distribution. In nonlinear filtering, this appears in the form of a
question about the stability of the filter~\cite{chigansky2009intrinsic, chigansky2006stability, crisan2011handbook}. A filter is stable if two
different initial distributions lead to the same filtering
distribution asymptotically in time. This is a desirable property for
a filter since the choice of the initial distribution is arbitrary and
we desire the filter to ``forget'' about this arbitrary choice. Thus
filter stability is an extensively studied topic, but mainly in the context of stochastic dynamics with only limited results for deterministic dynamical systems~{\color{mypink}\cite{reddy2021stability, oljavca2021exponential, reddy2019asymptotic, Cerou2000LongTB}}.

In practice, nonlinear filters need to be implemented numerically and we
will focus on two of the most commonly used methods, namely, particle
filters (PF) and ensemble Kalman filters (EnKF). The stability of
particle filters has been extensively studied~\cite{Chopin2020, whiteley2013stability} while very few theoretical results are known for the stability of EnKF~\cite{del2018stability}, though some results related to filter divergence (which is quite distinct from an unstable filter) and accuracy for the EnKF are available~\cite{KLS14, law2016filter, gottwald2013mechanism}. We
note that the assumptions used to prove stability of PF are
not satisfied by a deterministic dynamical model and thus their
stability in the context of data assimilation needs to be explored
numerically. This is the main aim of the present paper.

In order to assess filter stability directly using the definition
(see definition~\eqref{def-stab}), we need to compute distances between
probability distributions. This has been a challenging task, but recently proposed Sinkhorn algorithm provides an efficient method for this task~\cite{genevay2019entropy, genevay2018learning, feydy2019interpolating, arjovsky2017wasserstein}. One of the novelties of this paper is to demonstrate the use of the Sinkhorn algorithm in the context of data assimilation for directly studying stability.

Numerical studies of filter divergence, especially in the context of twin experiments where synthetic observations are generated using the model, have focused on assessing whether the filter remains bounded or whether the `error' or `bias' of the filter (commonly called RMSE, i.e. distance of the filter mean from the numerical trajectory which is used to generate the synthetic observations) remains bounded in time~\cite{KLS14, law2016filter}. But this does not provide a direct indication of stability as defined in~\eqref{def-stab}. In this paper, we demonstrate that there is a linear relationship between the filter error and the distance between two filters started with two different initial distributions, the latter giving a more appropriate measure of stability. This direct relation between filter stability and filter RMSE for particle and ensemble Kalman filters for deterministic dynamics is the other main contribution of this paper.

The outline of the paper is as follows. The next section~\ref{sec-problem} introduces the mathematical setting for the problem of filter stability. This is followed, in section~\ref{sec-method}, by a description of the dynamical model we use (Lorenz-96), the filtering algorithms (particle filter and ensemble Kalman filter), and the Sinkhorn algorithm for computing distances between probability distributions. The following section~\ref{sec-results} presents the numerical results, followed, in section~\ref{sec-conclude}, by a discussion including avenues for further research.

%% file: 02-problem.tex
\subsection{The nonlinear filtering problem}
In this paper we work with dynamical model given by a deterministic and chaotic ODE. Our model state is $d$-dimensional and the flow corresponding to the model is denoted by $\phi:\mathbb R \times \mathbb R^d\to\mathbb R^d$. We observe the model every $g$ units of time. So the model state $x_k$ follows a discrete-time deterministic dynamical system $f_g \stackrel{\rm def}{=} \phi(g, \cdot): \mathbb{R}^d \to \mathbb{R}^d$. The observations $y_k \in \mathbb{R}^q$ is related to the model state by the observation operator (a linear projection throughout this paper) $H: \mathbb{R}^d \to \mathbb{R}^q$ for $k = 0, 1, \dots$, as follows:
\begin{align}
&x_{k+1} = f_g(x_k), \quad x_0 \sim \mu \,, \label{eq-state}\\
&y_{k} = Hx_k + \eta_k \,, \quad  \label{eq-obs}
\end{align}
where $\mu$ is the initial distribution of the model state $x_0$ at time 0,  and $\eta_k \sim \mathcal{N}(0_q, \sigma^2I_q)$ are \emph{iid} Gaussian errors in the observation, and are assumed independent of $\mu$. Given observations $y_0, y_1, \cdots, y_n$, the goal of filtering is to estimate the conditional distribution of the model state at time $n$ conditioned on observations up to that time: $x_n|y_{0:n} \sim \pi_{n}(\mu)$, where the dependence on the distribution $\mu$ of the initial condition $x_0$ is made explicit since our focus will be on filter stability. %

\subsection{Filter stability}\label{ssec-def-stab} 
In practice we often do not know the initial distribution $\mu$. In such a case, when a different initial condition $\nu$ is chosen, one obtains a different filter, denoted by $\pi_{n}(\nu)$, by using the same set of observations and using the same algorithm. A measure of robustness of a filtering algorithm is how well it is able to "forget" the initial distribution, which motivates the following definitions.
{\color{mypink}

There are two different kinds of randomness that one needs to deal with in the setup above, the initial condition and the observation noise. Suppose $x_0:\Omega \to\mathbb R^d$ is our random initial condition. Consider $x_0(\omega)$, a realization of this initial condition. Now that we have fixed a realization of $x_0$, $\pi_n(\nu)$ and $\pi_n(\mu)$ become random measures whose randomness is determined only by the observation noise. For $x_0(\omega)$ we can compute the following expectation with respect to the observation noise 
\begin{align}
    \mathbb E \left|\int_{\mathbb R^d}  h(x)\, {{\pi_n(\mu, dx)}} - \int_{\mathbb R^d}  h(x)\, {{\pi_n(\nu, dx)}}\right|
\end{align}
for a bounded, continuous function $h$. If this expectation approaches $0$ as $n\to\infty$ for any such $h$ we can say that the filter is "pointwise" stable for the initial realization $x_0(\omega)$. And if the filter is "pointwise" stable for almost all realizations of $x_0$, we call the filter stable.  

\cite{reddy2021stability} explores the "true" filter stability for deterministic dynamics. Below we adapt the definition for numerical filter stability. Below $\hat\pi$ denotes the numerical approximation of the true filter $\pi$.
\begin{defn}[Stability-RA \cite{reddy2021stability}] A numerical filter is stable if  for any
measure $\nu$ with $\mu\ll\nu$ we have, 
\begin{align}
    \lim_{n\to\infty}\mathbb E \left|\int_{\mathbb R^d}  h(x)\, {{\hat\pi_n(\mu, dx)}} - \int_{\mathbb R^d}  h(x)\, {{\hat\pi_n(\nu, dx)}}\right|\label{def-stab-sumith} = 0 \,,
\end{align}
for any bounded and continuous $h$, $\mu$-almost everywhere in the sense described above.
\label{def-stab-ra} \end{defn}
Note that the expectation above is taken with respect to the observation noise only. Although \eqref{def-stab-sumith} captures the notion of filter stability quite well, from a computational perspective we can improve on it in the following aspects. 
\begin{itemize}
 \item Computing the expectation for every possible bounded and continuous function is infeasible.
\item In real world applications we might not have access to $\mu$ and therefore an expectation independent of $\mu$ is preferable.

\end{itemize}
In order to overcome above difficulties and to assess filter stability numerically, we devise the following definition which can be proven to be a stronger version of definition~\ref{def-stab-ra} in an appropriate sense (see theorem~\ref{thm:strength} in appendix).
\begin{defn}[Stability-MRA \cite{mandal2021stability}]A numerical filter is said to be stable if for any two distributions {$\nu_1, \nu_2$}, the following holds, 
\begin{align}
    {\lim_{n\to\infty}\mathbb E[D(\hat\pi_n(\nu_1), \hat\pi_n(\nu_2))]} = 0 \,,
\label{eq-stablaw} \end{align}
$\mu$-almost everywhere, where $D$ is a distance on $\mathcal P(\mathbb R^d)$, the space of probability measures on $\mathbb R^d$.
\label{def-stab} \end{defn}

Note that even with the modifications, definition~\ref{def-stab} remains hard to compute in the following aspects.
\begin{itemize}
    \item Computing the limit for every possible pair $\nu_1, \nu_2$ is infeasible.
    \item Computing the limit for every possible initial realization $x_0(\omega)$ is infeasible.
\end{itemize}

The last two difficulties also arise in definition~\ref{def-stab-ra} and are unavoidable in some sense in a complete definition of filter stability. But even with these difficulties we can explore numerical filter stability in a meaningful albeit slightly limited way.  Although we demonstrate results for a single realization $x_0(\omega)$ here, this realization was generated randomly and different initial realizations yield qualitatively similar results which are consistent with the stability definition~\ref{def-stab} and hence are not included in the paper to avoid repetition.
}

The main aim of this paper is to study the stability of two popular filtering algorithms, namely the particle filter (PF) and the ensemble Kalman filter (EnKF) by studying the limit in~\eqref{eq-stablaw}, where we choose the Wasserstein metric $W_2$ as our distance $D$ on $\mathcal P(\mathbb R^d)$. Previous work \cite{mandal2021stability} has shown these filters to be stable. Here we study the rate of convergence of the expectation in \eqref{eq-stablaw} and how it varies with respect to the time between the observations denoted by $g$ and the observational uncertainty or the error variance $\sigma^2$. Thus we study $\mathbb{E}[D(\hat{\pi}_n(\nu_1), \hat{\pi}_n(\nu_2))]$ as a function of time $n$ for PF and EnKF algorithms. {\color{mypink}In the following discussion we sometimes abuse the notation and use $\pi$ to mean $\hat \pi^{\rm PF}$ or $\hat \pi^{\rm EnKF}$ with clear context.} We now describe these numerical filtering algorithms, followed in section~\ref{ssec-sink} by a description of the Sinkhorn algorithm for computing distances between probability distributions.

\subsection{Ensemble Kalman Filters}\label{ssec-enkf}
Kalman filters provide the closed form solutions to the Bayesian filtering equations in the scenario when the dynamic and measurement models are linear Gaussian or if the state equation~\eqref{eq-state} looks like 
\begin{align}
   x_{k+1}=A_k x_k + \alpha_k \,, \label{eq-state-linear} 
\end{align}
where $\alpha_k\sim\mathcal N(\mathbf{0},Q_k)$. The filtering distribution in this special case turns out to be Gaussian. The mean and covariance of this distribution is computed recursively in two steps, a prediction step where the effect of the hidden dynamics is captured and $p(x_k|y_{1:k-1})$ is computed and an update step where the observation $y_k$ is taken into account to give the filtering distribution $p(x_k|y_{1:k})$ using Baye's rule and well-known properties of the multivariate Gaussian distribution.

Ensemble Kalman filters can be thought of as an approximation of the original Kalman filter where the filtering distribution is represented by a collection of particles, as is the norm in Monte Carlo-based methods. The ensemble representation is akin to  dimension reduction which leads to computational feasibility for systems with large state space dimension $d$.\cite{Evensen03}. Localization, which is the process of weeding out long range spurious correlations, has made EnKF more applicable as well as wildly popular in high-dimensional data assimilation problems for spatially extended systems. For a discussion about localization see \cite{carrassi2018data}. We use Gaspari-Cohn function as our choice of localization function {\color{mypink}with radius set to 2}. 
\begin{algorithm}[!t]
\textcolor{mypink}{Initialize $N$ particles $\{x_0^i\}_{i=1}^N$ according to the initial distribution and set $x_0^{i,a}=x_0^i$ \\
Set $\rho$ as the Gaspari-Cohn localization matrix \cite{carrassi2018data}.\\
\For{$k=1,\cdots,n$}{
    \For{$i=1,\cdots,N$}{
        $x^{i,f}_{k}\leftarrow f_g(x_{k-1}^{i,a})$
    }
    $m_{k}^{f}\leftarrow \frac{1}{N}\sum_{i}x_{k}^{i,f}$\\ %
    $P_{k}^{f}\leftarrow \rho \circ \frac{\sum_{i}\left(x_{k}^{i,f}-m_{k}^{f}\right)\left(x_{k}^{i,f}-m_{k}^{f}\right)^\top}{N-1}$\\
    $K \leftarrow P^{f}_{k}H^{T} \left[HP^{f}_{k} H^{T}+R_{k}\right]^{-1}$\\
    \For{$i=1,\cdots,N$}{
        Sample $\eta^{i}_{k} \sim \mathcal{N}(0_q, \sigma^2I_q)$\\
        $y^{i}_{k}\leftarrow y_{k}+\eta^{i}_{k}$\\
        $x^{i,a}_{k} \leftarrow x^{i,f}_{k}+K\left[y^{i}_{k}-Hx^{i,f}_{k} \right]$
        }
    $\hat\pi_k \leftarrow\frac{1}{N}\sum_{i=1}^N \delta_{x^{i, a}_k}$
}
\caption{EnKF with covariance localization in state-space. $\circ$ denotes Hadamard product.}
\label{algo-enkf} 
}
\end{algorithm}
{\color{mypink} The details of the exact implementation that we use can be found in algorithm~\ref{algo-enkf}.}

\subsection{Particle Filters}\label{ssec-pf}
Particle-filters are also Monte Carlo-based filters that recursively compute importance sampling approximations of the filtering distribution $p(x_k|y_{1:k})$. PFs also follow the Bayesian paradigm of two-step recursion with prediction and update steps. The filtering distribution is represented as a collection of weighted particles. In the prediction step the particles are evolved in time according to $\eqref{eq-state}$ which gives us the prior for the next Bayesian update step where the weights are adjusted appropriately to account for the observation. For an excellent overview of the PF algorithm see \cite{doucet2009tutorial}. PFs do not rely on linearity or Gaussianity of dynamic of observation models which make them powerful but unless the number of particles scale exponentially with $d$, PFs experience weight degeneracy and provide poor estimates \cite{bengtsson1981dynamic}. In order to combat weight degeneracy, a resampling step is performed after the Bayesian update where particles with negligible weights are replaced with particles with higher weights. Many variants of the standard or the bootstrap PF have been proposed and the interested reader can see \cite{farchi2018comparison} for a discussion. However, applying PFs on problems with significantly high dimensions still remains a challenge. 

We use the bootstrap particle filter for our experiments with a stochastic resampling step where we place a Gaussian distribution with a pre-determined, small covariance {\color{mypink}$\tilde\sigma^2$ (set to be $0.5$ in the actual experiments)} around the best-performing particles and sample new particles according to the weights. 
{\color{mypink} The details of the exact algorithm can be found in algorithm~\ref{algo-bpf}.}
\begin{algorithm}[!t]
\textcolor{mypink}{Initialize $N$ particles $\{x_0^i\}_{i=1}^N$ according to the initial distribution with equal weights $\left\{w_0^i=\frac{1}{N}\right\}_{i=1}^N$. Set $\tilde\sigma$. Below $S[i]$ denotes $i$-th element of $S$.\\
 \For{$k=0,\cdots,n$}{
      \If{$k>0$}{
      \For{$i=1,\cdots,N$}{
            $x^i_k\leftarrow f_g(S[i])$
            }
      }
      Sample $u\sim\mathcal{U}\left(0, \frac{1}{N}\right)$\\
      \For{$i=1,\cdots,N$}{
            $w^i_k\leftarrow p(y_k|x_k^i)$\\
            $U_{i} \leftarrow u + \frac{i-1}{N}$
            }
      $W\leftarrow\sum_{i=1}^Nw^i_k$\\
      \For{$i=1,\cdots,N$}{ 
            \If{$|\{U_j:\sum_{l=1}^{i-1}w^l_k \le WU_j\le\sum_{l=1}^{i}w^l_k\}|>0$}{
                    tag $x^i_k$ as significant
            }
            }
      Set $S\leftarrow\{x^{i_1}_k, x^{i_2}_k, \cdots, x^{i_m}_k\}$ as the set of  significant particles and compute $N_j\propto w^{i_j}_k:  \sum_{j=1}^mN_j = N$.\\
      \For{$j=1, \cdots, m$}{
            $S\leftarrow S\cup\{N_j-1\text{ samples from }\mathcal{N}(x^{i_j}_k, \tilde\sigma^2I_d)\}$
            }
    $\hat\pi_k \leftarrow\frac{1}{N}\sum_{i=1}^N\delta_{S[i]}$
  }
 \caption{BPF with offspring-based resampling}\label{algo-bpf}}
\end{algorithm}

{\color{mypink}\subsection{Choice of distance $D$} Although the stability definitions are independent of  any kind of specific distance and can be computed with other choices of $D$ we choose $D$ to be the $2$nd Wasserstein distance $W_2$. We justify our choice with the reasons below.
\begin{itemize}
    \item The is an efficient algorithm for approximating $W_p$ which is described below. 
    \item Some nice geometric properties of $W_p$ e.g. metrizing the convergence in law \cite{feydy2019interpolating} are missing from other distances or distance-substitutes like the  total variation distance and  the KL-divergence. 
    \item Moreover, $W_p$ does not require the notion of absolute continuity unlike KL-divergence or Hellinger distance which is useful for comparing empirical distributions.
\end{itemize}
For a comparison of these distances the interested reader may see \cite{arjovsky2017wasserstein} where example $1$ (learning parallel lines) depicts how the output of $W_p$ can often be intuitive because Wasserstein distances lift the standard metrics on $\mathbb R^d$ to the probability space $\mathcal P(\mathbb R^d)$ unlike KL-divergence or total variation. Lastly, we use $p=2$ for no reason other than the familiarity of the $2$-norm on Euclidean spaces.}

\subsection{Sinkhorn divergence} \label{ssec-sink}

The $p$-th Wasserstein distance ($W_p$) between probability measures with $p$-th finite moment on metric spaces have many desirable geometric features which stems form the fact that its definition extends the distance function on the metric space to a distance on the space of probability measures on the metric space. For a discussion see \cite{feydy2019interpolating, arjovsky2017wasserstein}. $W_1$ or the earth mover's distance has been used in various problems e.g. comparing colour histograms, solving resource allocation problems etc. When applied to two sampling distributions with both having sample size $k$, computing $W_1$ is equivalent to solving a constrained linear programming problem in $n = k^2$ variables. Since LPPs take $O(n^3)$ time to solve a problem with $n$ variables, computing $W_1$ takes $O(k^6)$ time which is prohibitively expensive.

In recent years it has been noted that by regularizing the optimization problem that defines the Wasserstein distance, one can attempt to solve the dual to the regularized problem which is akin to solving a convex optimization problem. For a comprehensive discussion see \cite{genevay2019entropy}. The dual problem can be solved using a variant of the Sinkhorn-Knopp algorithm for finding a doubly-stochastic matrix given a square matrix with positive entries. The solution to the regularized problem is known as the Sinkhorn divergence since it fails to satisfy the triangle inequality and is not an exact distance on the space of probability measures.

{\color{mypink}Here we focus on the case $p=2$.} For two probability measures $\mu$ and $\nu$ on $\mathbb R^d$ with finite first and second moments, the Sinkhorn divergence $S_\varepsilon$ is defined as follows \cite{feydy2019interpolating}.
\begin{align}
    &\text{OT}_\varepsilon(\mu, \nu) \stackrel{\text{def}}{=} \min_{\pi \in \mathbb{S}}\left[\int\|x-y\|_2^2\,d\pi(x, y) + \varepsilon\text{KL}(\pi|\mu\otimes\nu)\right] \,, \label{def-ot}\\
    &\text{S}_\varepsilon(\mu, \nu) \stackrel{\text{def}}{=} \text{OT}_\varepsilon(\mu, \nu) -\frac{1}{2}\text{OT}_\varepsilon(\mu, \mu)-\frac{1}{2}\text{OT}_\varepsilon(\nu, \nu) \,, \label{def-sink}
\end{align}
where the minimisation is over the set $\mathbb{S}$ of distributions $\pi$ with the first and second marginals being $\mu$ and $\nu$ respectively  
and $\text{KL}$ is the Kullback–Leibler divergence.
Moreover, it turns out \cite{feydy2019interpolating} that
\begin{align}
    \lim_{\varepsilon \to0}\sqrt{S_\varepsilon(\mu, \nu)} = W_2(\mu, \nu) \,,\label{eq-sinkhorn-limit}
\end{align}
and therefore for small enough $\varepsilon$ we 
obtain a good approximation of $W_2$. We use the following notation for this approximation: $D_\varepsilon = \sqrt{S_\varepsilon}$.

In our experiments we compute $S_\varepsilon(\mu, \nu)$ for sampling distributions $\mu=\sum_{i=1}^N\mu_i\delta_{x_i}$ and $\nu=\sum_{j=1}^M\nu_j\delta_{y_j}$ with $\varepsilon=0.01$ where $\{x_i\}_{i=1}^N$ and $\{y_j\}_{j=1}^M$ are points in $\mathbb R^d$. {\color{mypink} A detailed justification of this choice of $\varepsilon$ is in appendix~\ref{ssec-dw}.} The Sinkhorn divergence algorithm being a fixed point iteration is extremely fast. The exact procedure is given in algorithm~\ref{algo-sink}. The authors of \cite{feydy2019interpolating} show that the algorithm is parallelizable with respect to sample-size. The dimension dependence of the algorithm is only explicitly apparent while calculating the distance matrix and consequently the algorithm itself scales only linearly in $d$. But to accurately represent the underlying distribution with increasing dimension one would require to compute $S_\varepsilon$ with increasing sample size. For a detailed discussion of sample complexity of the Sinkhorn divergence see chapter 3 of \cite{genevay2019entropy}.
\begin{algorithm}[!t]
 \hspace*{\algorithmicindent} \textbf{Input: } $\{\mu_i\}_{i=1}^N, \{x_i\}_{i=1}^N, \{\nu_j\}_{j=1}^M, \{y_j\}_{j=1}^M$ \\
 \hspace*{\algorithmicindent} \textbf{Output: } $S_\varepsilon\left(\sum_{i=1}^N\mu_i\delta_{x_i},\sum_{j=1}^M\nu_j\delta_{y_j} \right)$\\ 
Note the definition, $\text{LSE}_{k=1}^LV_k\stackrel{\text{def}}{=}\log\sum_{k=1}^L\exp(V_k)$.\\
Initialize $a_i\leftarrow0\;\forall\; i=1,\cdots,N$ and $ b_j\leftarrow 0,\;\forall\; j=1,\cdots,M$.\\
iteration $\leftarrow 0$\\
\While{$\min$\{$L_1$ relative errors in $a$ and $b$\} $> 0.1\%$  }{
    \For{$i=1,\cdots, N$}{
        $a_i \leftarrow\
        -\varepsilon\text{LSE}_{k=1}^M\left(\log\nu_k+\frac{1}{\varepsilon}b_k - \frac{1}{\varepsilon}\|x_i-y_k\|_2^2 \right)$} %
    \For{$j=1,\cdots, M$}{
        $b_j \leftarrow\
        -\varepsilon\text{LSE}_{k=1}^N\left(\log\mu_k+\frac{1}{\varepsilon}a_k - \frac{1}{\varepsilon}\|x_k-y_j\|_2^2 \right)$}\
iteration $\leftarrow$   iteration + 1}\
$\text{OT}_{\mu, \nu}\leftarrow\sum_{i=1}^N\mu_i a_i +\sum_{j=1}^M \nu_j b_j$\\
Initialize $a_i\leftarrow0\;\forall\; i=1,\cdots,N$ and $ b_j\leftarrow 0,\;\forall\; j=1,\cdots,M$.\\
\While{$L_1$ relative error in $a > 0.1\%$}{
    \For{$i=1,\cdots, N$}{
        $a_i \leftarrow
        \frac{1}{2}\left[a_i - \varepsilon\text{LSE}_{k=1}^N\left(\log\mu_k+\frac{1}{\varepsilon}a_k - \frac{1}{\varepsilon}\|x_i-x_k\|_2^2 \right)\right]$}
    } %
\While{$L_1$ relative error in $b > 0.1\%$}{
    \For{$j=1,\cdots, M$}{
        $b_j \leftarrow\
        \frac{1}{2}\left[b_j - \varepsilon\text{LSE}_{k=1}^M\left(\log\nu_k+\frac{1}{\varepsilon}b_k - \frac{1}{\varepsilon}\|y_j-y_k\|_2^2 \right)\right]$}
    } %
$S_\varepsilon\leftarrow \text{OT}_{\mu, \nu} - \sum_{i=1}^N\mu_i a_i -\sum_{j=1}^M \nu_j b_j$
 \caption{Computation of $S_\varepsilon$}
\label{algo-sink}
\end{algorithm}

%% file: 03-method.tex
\subsection{Model} \label{ssec-models}
Proposed first in 1995 by Edward N. Lorenz, the Lorenz-96 equations are a set of autonomous equations said to be mimicking the circulation of the earth's atmosphere in an over-simplified manner. Although simple, they have had significant impact on the development of the dynamical systems theory, especially because of their chaotic nature in arbitrary dimensions. Since an important application of data assimilation is numerical weather prediction, the Lorenz systems are a natural first choice for experiments. Since their inception they have been extensively used in data assimilation literature. Here we use $d = 10$~dimensional Lorenz-96~\cite{Lorenz96, kekem2018dynamics} with forcing constant $F=10$. We observe the system $g$ units of time which fixes the evolution function $f$. We observe alternate coordinates starting from the first coordinate, so
\begin{align}
    y_{k,j} = x_{k,2j-1} + \eta_{k,j} \,,
\label{eq-altobs} \end{align}
for $j=1, 2, \cdots, q = \left\lceil\frac{d}{2}\right\rceil$ and $\eta_{k, j}\sim\mathcal N(0, \sigma^2)$. Throughout the paper, we use $\sigma^2=0.2, 0.4, 0.8, 1.6$ and \newline $g=0.01, 0.03, 0.05, 0.07, 0.09$. The choice of the dimension $d=10$  makes sure that we get reasonable performances from both the EnKF and the particle filter.

\subsection{Data generation}
Lorenz systems are known to have attracting sets. In  this paper, we focus on the special case when the filtering distributions are expected to be supported on the attractor. Not only does it mimic real world scenarios, it also lets us make use of the theory optimal transport distances, outlined in \cite{feydy2019interpolating}, when the probability measures are supported on a compact domain, since the Lorenz attractors are bounded sets. So we begin by finding a point on the attractor by randomly generating an initial point and evolving it according to $f_g$ for $10^5$ iterations. Starting from this point $x_0^{\text{true}}$ on the attractor, we generate a true trajectory according to \eqref{eq-state} and then generate $10$ different observation realizations for the same trajectory according to \eqref{eq-altobs} in order to compute the expectation over observational noise, as in \eqref{eq-stablaw}. {\color{mypink} For a justification of why $10$ observation realizations suffice for our study, see appendix~\ref{ssec-sample-size}}. 

\subsection{Initial distributions}\label{ssec-init-dist}
We use two Gaussian initial conditions. The first one $\mu_0$ is centered at the true state with a small variance representing the case when our guess for the initial distribution is unbiased and precise. Thus we expect the filter to continue to have those properties {\color{mypink} upto some time.} The second one $\mu_b$ is centered away from the true state with a significantly larger variance representing the case when our guess for the initial distribution is biased and imprecise. They are given by,
\begin{align}
    &\mu_0 = \mathcal{N}(x_0^{\text{true}}, 0.1\times I_d) \,, \nonumber \\
    &\mu_b = \mathcal{N}(x_0^{\text{true}} + 4\times1_d, I_d) \,,
\label{eq-3ic} \end{align}
where $1_d$ is a $d$-dimensional vector with all entries $1$. {\color{mypink}With this notation $x_0^{\rm true}$ corresponds to $x_0(\omega)$ in subsection \ref{ssec-def-stab}. Note that different realizations of $x_0$ produce similar results as shown here.}

\subsection{Metrics for filter stability}

To probe filter stability directly using the definition~\ref{def-stab}, we study the Sinkhorn divergence $\mathbb E[D_\varepsilon(\pi_n(\mu_0), \pi_n(\mu_b))]$ as a function of time. It has been well-known that in nonlinear filtering problems where the dynamic model is stochastic, under suitable additional conditions, the filter is exponentially stable and from an incorrect initial condition, it reaches stability in an exponential fashion (see, e.g., chapter 3 of~\cite{van2008hidden}). Although such results are not available for the case of deterministic dynamics, exponential decay is a natural or at least desirable behaviour for the temporal behaviour of the distance between two filters starting from different initial distributions. To explore this qualitatively, we fit a curve of the following form
\begin{align}
    \mathbb E[D_\varepsilon(\pi_n(\mu_0), \pi_n(\mu_b))] = a\exp(-\lambda t) + c \,, \label{eq:fit}
\end{align}
where time $t=$ assimilation step $\times$ observation gap = $ng$. One of the motivation is to understand whether the exponent $\lambda$ is related the dynamical quantities such as the Lyapunov exponents of the chaotic dynamical system under consideration.

In addition to stability, we also explore its relationship to the convergence of the filter mean toward the true signal as well as the uncertainty of the mean estimate. Motivated by the results about bounds on the former of these two in~\cite[Theorem~4.4]{KLS14} and \cite[Theorem~4.6]{law2016filter}, we define the following two quantities.

The first quantity aims to capture the bias of the filter and is the scaled $l_2$ error denoted by $e_n(\nu)$.
\begin{align}
    e_n(\nu) &\stackrel{\rm def}{=} \frac{1}{\sqrt{d}} \left\| \mathbb{E}_{\hat\pi_n(\nu)}[x_n] - \phi\left(ng, x_0^{\text{true}}\right) \right\|_2 \,, \nonumber \\
    &= \left[ \frac{1}{d} \sum_{i=1}^d \left(\frac{1}{N} \sum_{\alpha=1}^N x_n^{\alpha,i} - x_n^{\text{true}, i} \right)^2 \right]^{1/2} \,, \label{eq-error}
\end{align}
where $x_n^{\alpha,i}$ denotes the $i$-th coordinate of the $\alpha$-th member of the ensemble representing the filtering distribution at time $n$. Thus, $e_n(\nu)$ is the distance between the true state and the filter mean divided by square root of the state space dimension $d$, with $n$ denoting the assimilation step and $\nu$ the initial distribution of the filter. Note that from the results \cite[Theorem~4.4]{KLS14} and \cite[Theorem~4.6]{law2016filter} mentioned earlier, we expect $\mathbb{E}[e_n^2(\nu)] \sim \sigma^2$ asymptotically in time.

The second quantity $s_n(\nu)$ captures the uncertainty of the filter estimate.
\begin{align}
    s_n(\nu) & \stackrel{\rm def}{=} \left[ \frac{1}{d} \tr \left[ \mathbb{E}_{\hat\pi_n(\nu)}[\left(x_n - \mathbb{E}_{\hat\pi_n(\nu)}[x_n] \right)\left(x_n - \mathbb{E}_{\hat\pi_n(\nu)}[x_n] \right)^t \right]^{1/2} \right] \,, \nonumber \\ 
    &= \left[ \frac{1}{d} \sum_{i=1}^d \frac{1}{N-1} \sum_{\alpha=1}^N \left(x_n^{\alpha,i} - \sum_{\beta=1}^N x_n^{\beta,i} \right)^2 \right]^{1/2}\,, \label{eq-var}
\end{align}
 Thus, $s_n(\nu)$ is the square root of the trace of the sample covarianace of the filter. We note that we are not aware of any theoretical results that give any indication about the asymptotic in time limit of this quantity but it is reasonable to explore their relation with the observational uncertainty.

%% file: 045-filter-results.tex
We now present the numerical results on stability for PF and EnKF, as well as the dependence of the exponential rates on the observation gap $g$ and the observation noise strength $\sigma$. In the following figures~\ref{fig:bpf-enkf-fixed-ocov}-\ref{fig:bpf-enkf-fixed-ogap}, in top row in each of the figures, the dots represent distance $D_\varepsilon(\pi_n(\mu_0), \pi_n(\mu_b))$ between the posterior distributions at time $t=ng$ with the initial distribution $\mu_0$ and $\mu_b$ mentioned in~\eqref{eq-3ic} versus time for 10 realizations of the observational noise. We also plot the mean $\mathbb{E}[D_\varepsilon]$ averaged over these 10 observational realizations, and the best-fit curve~\eqref{eq:fit} in the same figures. Rows 2 and 3 contain, respectively, the expectation value of scaled $l_2$ error $\mathbb{E}[e_n^2]$ and uncertainty $\mathbb{E}[s_n^2]$ defined in~\eqref{eq-error} and \eqref{eq-var} versus time. Row 4 contains scatter plot of RMSE, defined as the square root of the expected value of $e_n^2$, averaged over the 10 observational realizations, for the filter with biased initial distribution $\mu_b$ versus the mean $D_\varepsilon$ between posteriors from the biased and the unbiased initial distributions.

\subsection{Dependence on observation gap}
\begin{figure}[t!]
\centering
Particle filter $\qquad \qquad \qquad \qquad \qquad \qquad \qquad \qquad $ EnKF\\
    \includegraphics[width=0.48\textwidth]{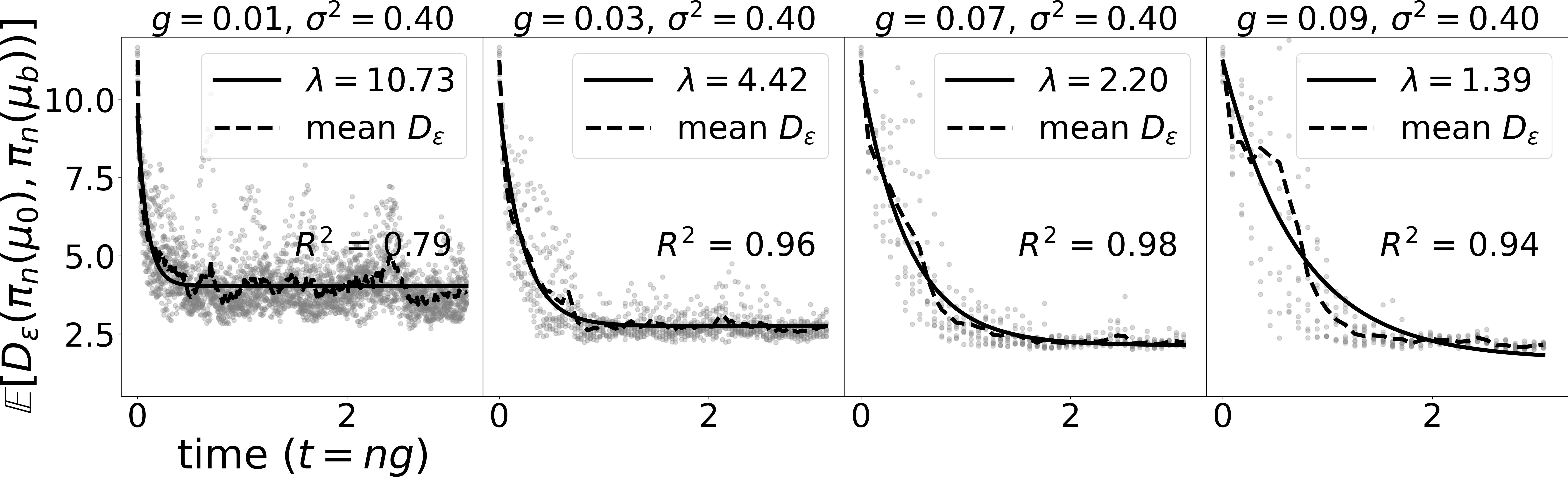} $\ \ $
    \includegraphics[width=0.48\textwidth]{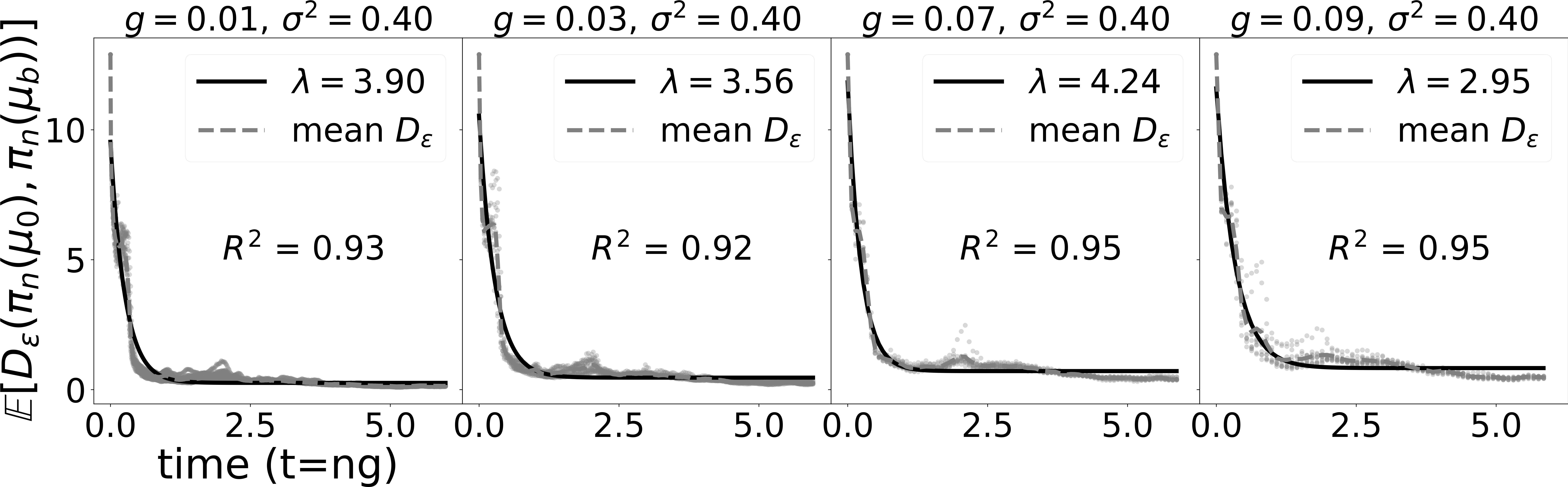}\\
   \includegraphics[width=0.48\columnwidth]{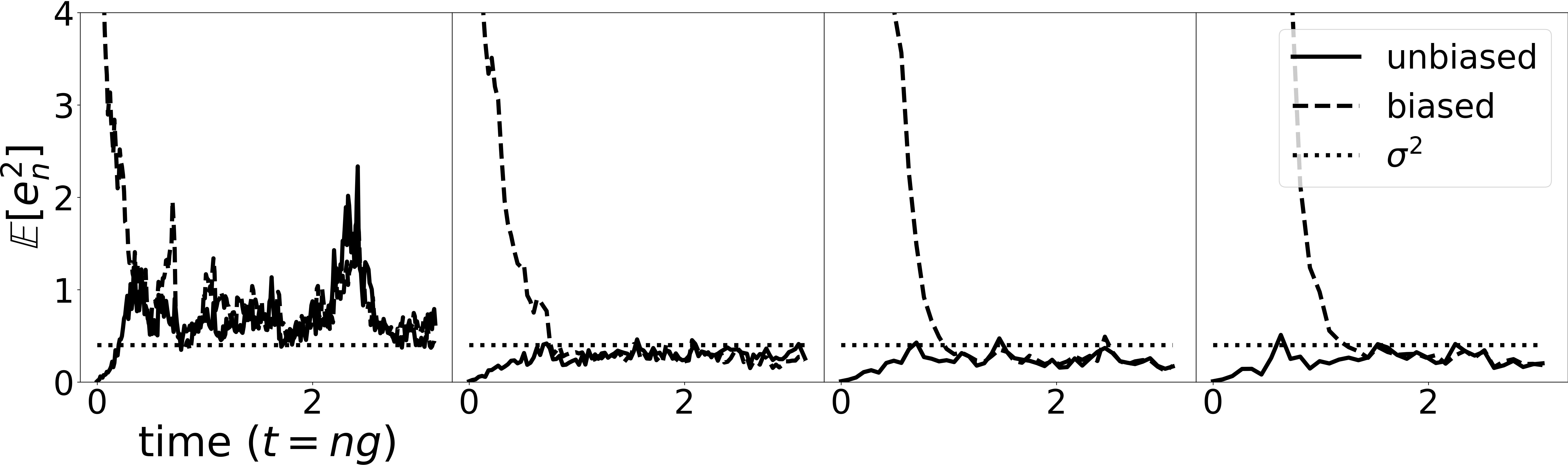} $\ \ $
    \includegraphics[width=0.48\textwidth]{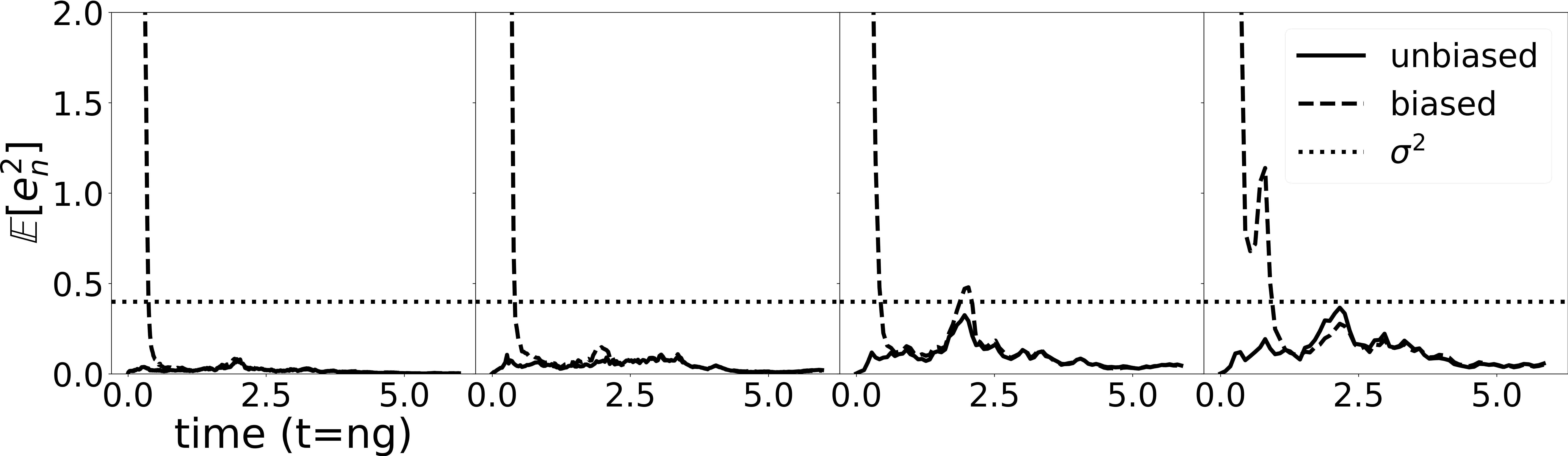}\\
   \includegraphics[width=0.48\columnwidth]{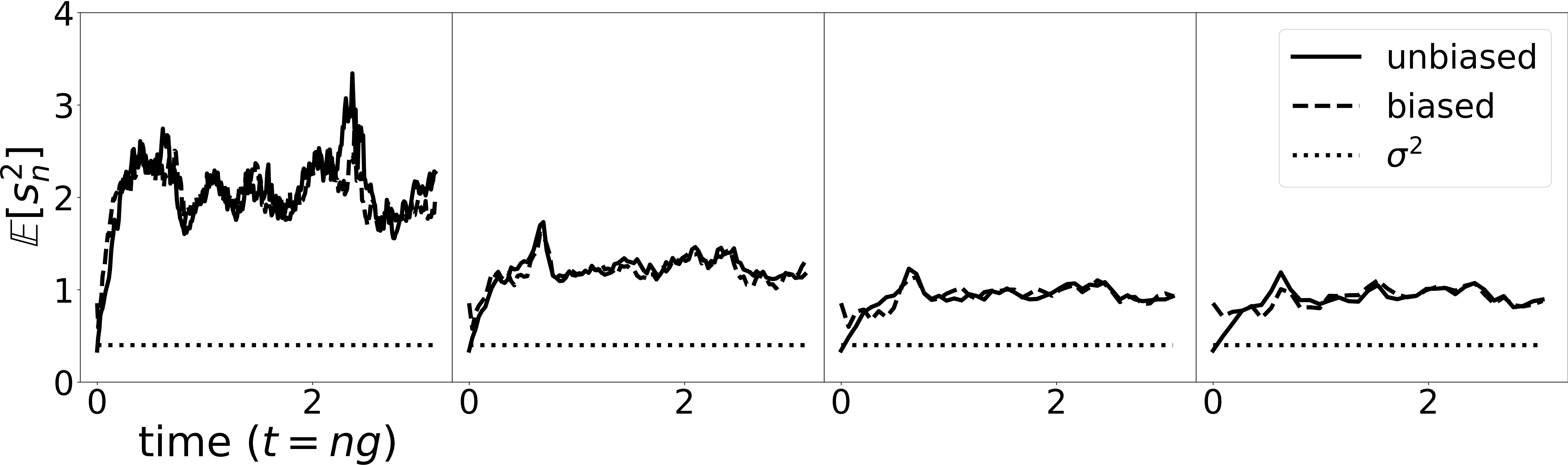} $\ \ $
    \includegraphics[width=0.48\textwidth]{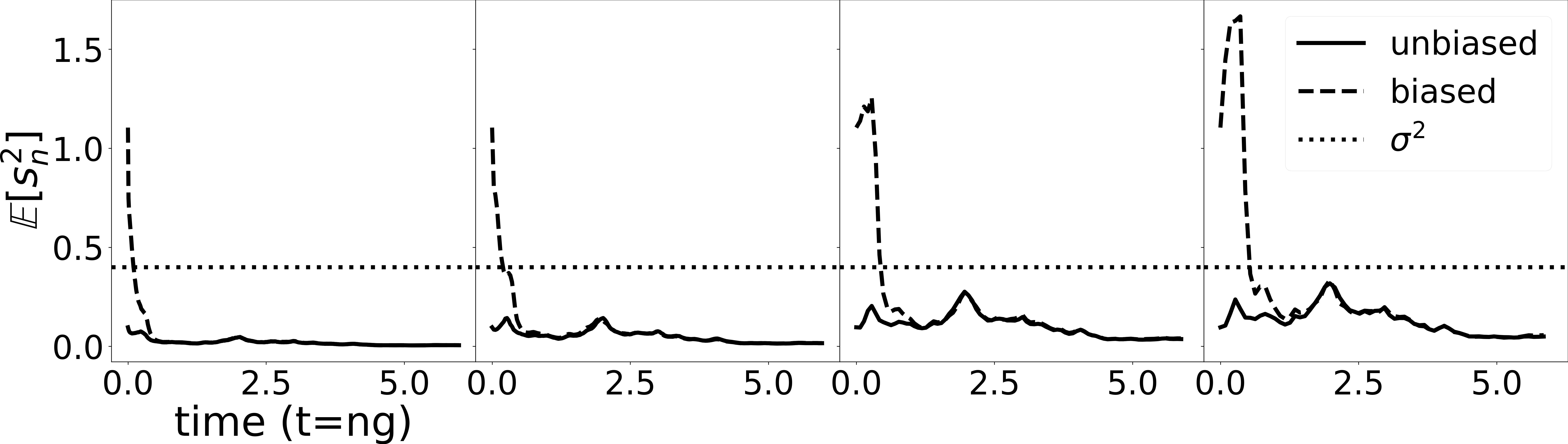}\\
   \includegraphics[width=0.48\columnwidth]{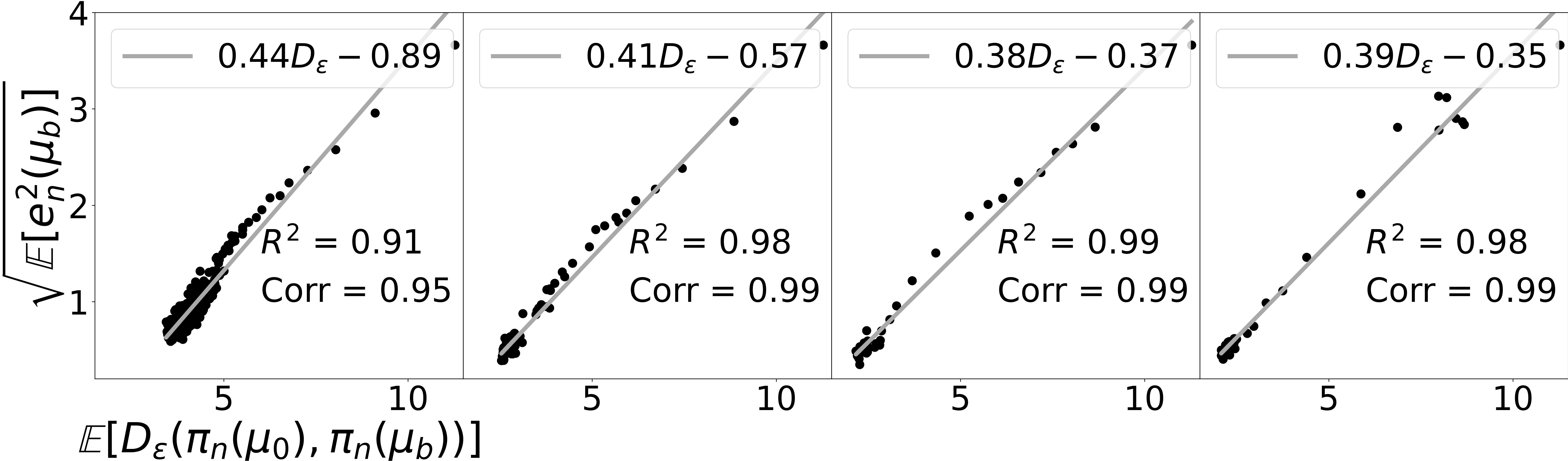} $\ \ $
    \includegraphics[width=0.48\textwidth]{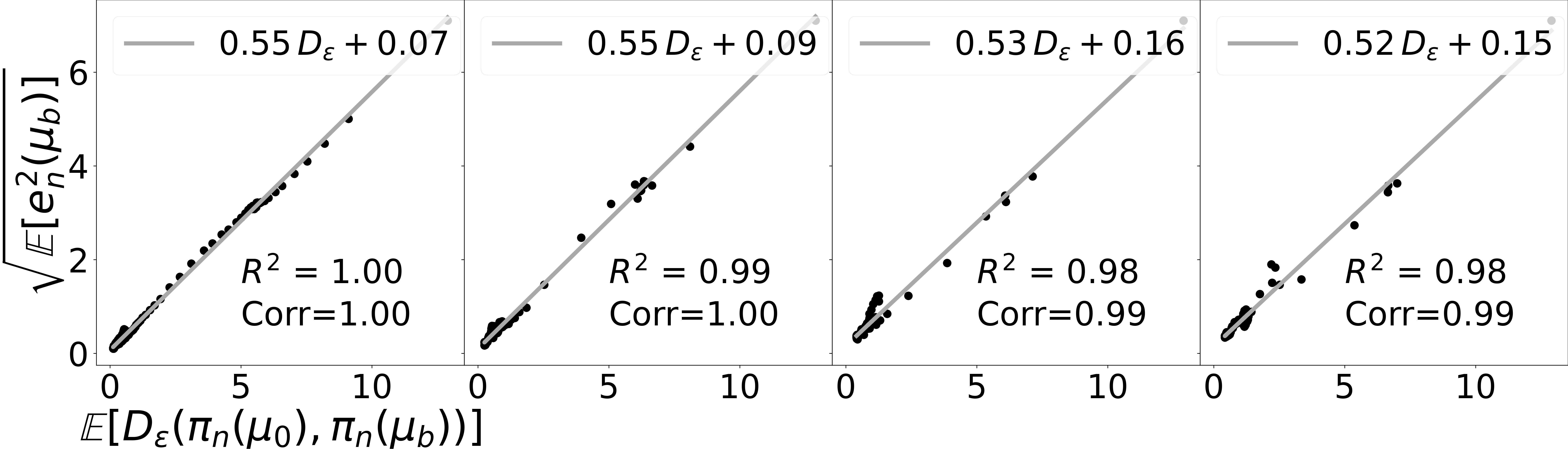}
\caption{The left and right panels show the results for PF and EnKF respectively with fixed observational error variance $\sigma^2 = 0.4$, and each column contains the results for different time between observations $g = 0.01, \, 0.03, \, 0.07, \, 0.09$.
Row 1: Mean $D_\varepsilon$ versus time. The dots represent 10 different realisations. The solid line is the exponential best-fit line for the mean $D_\varepsilon$ as in~\eqref{eq:fit}.
Row 2: Mean scaled $l_2$ error from~\eqref{eq-error} versus time for the two initial distributions.
Row 3: Mean uncertainty from~\eqref{eq-var} versus time for the two initial distributions. The constant dotted line in rows 2 and 3 shows the observational error variance $\sigma^2$ for reference.
Row 4: RMSE versus mean $D_\varepsilon$. Pearson correlation coefficient between these two quantities is depicted alongside the goodness of fit for the best-fit line.}
\label{fig:bpf-enkf-fixed-ocov}
\end{figure}

We first discuss the results of assimilation with both PF and EnKF for a fixed observation covariance $\sigma^2 = 0.4 I$ with varying observation gap. As shown in figure~\ref{fig:bpf-enkf-fixed-ocov}, the mean $D_\varepsilon$ falls exponentially over time until reaching a stationary value for both the filters. Table~\ref{table:fixcov} shows the values of the coefficients of the best-fit of mean $D_\varepsilon$ versus time, according to~\eqref{eq:fit}, for different observation gaps $g$.

\begin{table}[t!]
\centering
\begin{tabular}{|c|c|c|c|c|c|} 
 \hline
\multicolumn{2}{|c|}{$g$} & $\bm{0.01}$ & $ \bm{0.03}$  & $\bm{0.07} $ & $\bm{0.09}$ \\ [0.5ex] 
\hline
\multirow{2}{*}{\text{a}} & \textbf{PF}& 5.367 $\pm$ 0.080 & 7.077 $\pm$ 0.055 & 8.672 $\pm$ 0.084 & 9.54 $\pm$ 0.40 \\\cline{2-6}
& \textbf{EnKF}& 9.28 $\pm$ 0.13 & 10.08 $\pm$ 0.27 & 11.11 $\pm$ 0.32 & 10.76 $\pm$ 0.37 \\
\hline
\multirow{2}{*}{$\lambda$}& \textbf{PF} & 10.73 $\pm$ 0.76 &  4.423 $\pm$ 0.058 & 2.203 $\pm$ 0.021 & 1.392 $\pm$ 0.052 \\ \cline{2-6}
& \textbf{EnKF} & 3.904 $\pm$ 0.085 &  3.56 $\pm$ 0.15 & 4.24 $\pm$ 0.21 & 2.95 $\pm$ 0.17 \\
\hline
\multirow{2}{*}{c} & \textbf{PF} & 4.03425 $\pm$ 0.00083 & 2.7524 $\pm$ 0.0018 & 2.1362 $\pm$ 0.0093 & 1.69 $\pm$ 0.14\\ \cline{2-6}
& \textbf{EnKF} & 0.258 $\pm$ 0.016 & 0.459 $\pm$ 0.036 & 0.711 $\pm$ 0.046 & 0.827 $\pm$ 0.064\\
\hline
\end{tabular}
\caption{Parameters of the best-fit for the mean $D_\varepsilon$ versus time as in~\eqref{eq:fit} with associated confidence intervals for fixed observation covariance $\sigma^2 = 0.4$ and different observation gap $g$ shown in the top row.}
\label{table:fixcov}
\end{table}

For PF, we see that the rate $\lambda$ decreases with increasing observation gap $g$ but for EnKF, the rates are not significantly affected by the change in $g$. The highest Lyapunov exponent for the model (with the chosen parameter value) is approximately $\lambda_{\text{max}} = 1.7$ whereas the exponential rate $\lambda$ for EnKF is seen to be in the range of $(3.0, 4.2)$, close to $2 \lambda_{\text{max}}$, indicating a possible close relation between the dynamics and the EnKF that could be explored further. The exponential rate for the PF does not seem to show such a relation.

Another difference between PF and EnKF may also be noted: with increasing observation gap, the stationary value $c$ for the PF decreases whereas it increases for the EnKF. Also, the stationary values $c$ of the $D_\varepsilon$ for EnKF are significantly lower compared to corresponding values for PF. These asymptotic values of $D_\varepsilon$ over time for both the filters can be explained by their mean posterior covariance using the following argument.

A characteristic of the numerical distance $D_\varepsilon$ is that for two different {\it i.i.d.}~samples drawn from the same probability distribution, $D_\varepsilon$ has a nonzero positive value. Statistically, $D_\varepsilon$ between two empirical measures approach this value at which they are essentially representing the same distribution and cannot be distinguished. For a fixed dimension $d$ and sample size $N$, this value increases with increasing covariance of the distribution \cite[Figure~1 and discussion therein]{mandal2021stability}.

The mean posterior covariance trace is directly proportional to the $s_n^2$. With increasing observation gap $g$, the mean uncertainty decreases for PF while it increases for EnKF. Hence the asymptotic value of $D_\varepsilon$ decreases with increasing observation gap for PF, while for the latter, it increases. We note that the previous paragraph explains the asymptotic value of $D_\varepsilon$, but the difference in the behaviour of the filter uncertainty $s_n$ for PF and EnKF as a function of observation gap needs to be explored further.

The scaled $l_2$ errors also reach an asymptotically constant value around the same time the corresponding filters stabilize in $D_\varepsilon$. The scatter plots for the RMSE against the $D_\varepsilon$ shows strong correlation between them. This suggests that we can use the RMSE over time as a good indicator for the time when the filter stabilizes. Note that the methods in this paper give us a direct way to check whether a numerical filter is stable for a given dynamical and observational model, and the relation between the filter stability and the $l_2$ error or bias $e_n$ implies that a stable filter may be expected to be an accurate one.

We note that in the plots in the bottom row, the cluster at the bottom left corresponds to the time after which both RMSE and the $D_\varepsilon$ have reached their stationary values. Even for two different biased initial distributions for the filter, there is a finite transient growth after which the $D_\varepsilon$ falls exponentially \cite{mandal2021stability}. Although not shown here, the linear regime is still present in the scatter plot of the RMSE of either one of them versus the $D_\varepsilon$ in those cases. 

\subsection{Dependence on observation noise $\sigma$ }

\begin{figure}[t!]
\centering
Particle filter $\qquad \qquad \qquad \qquad \qquad \qquad \qquad \qquad $ EnKF\\
    \includegraphics[width=0.48\textwidth]{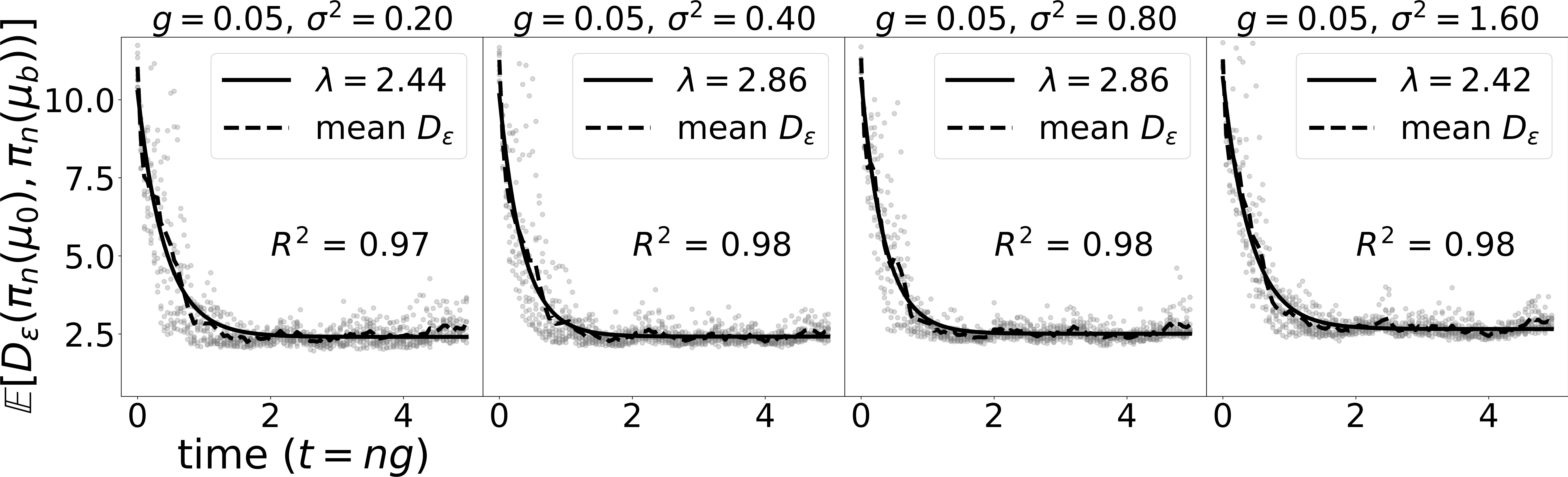} $\ \ $
    \includegraphics[width=0.48\textwidth]{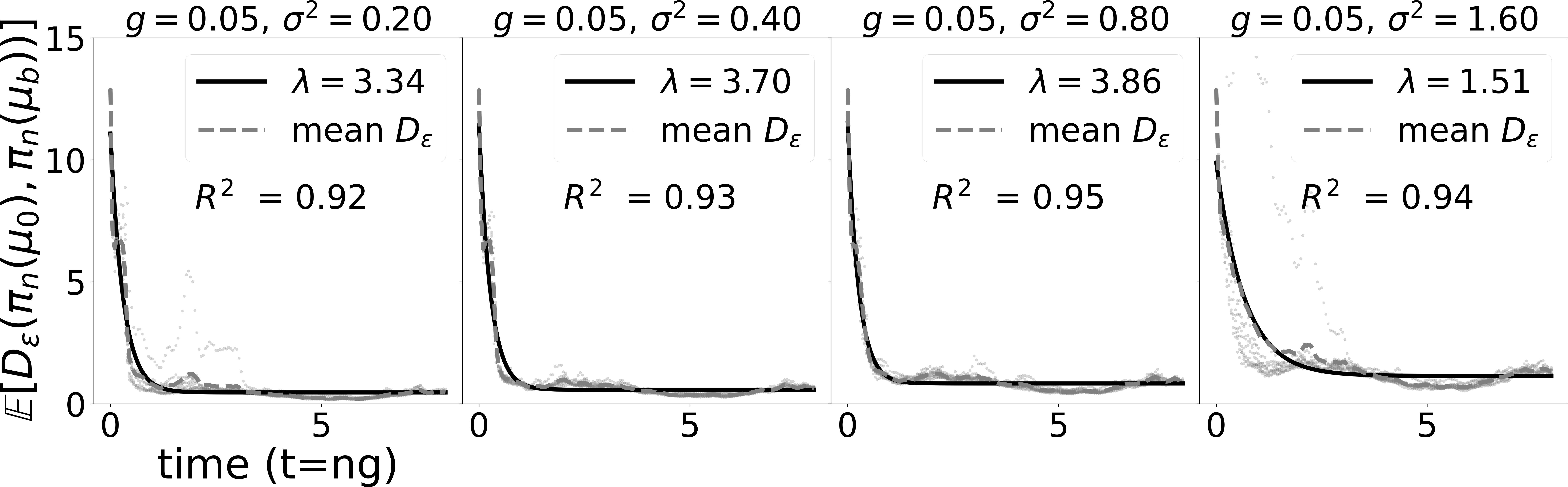}\\
   \includegraphics[width=0.48\columnwidth]{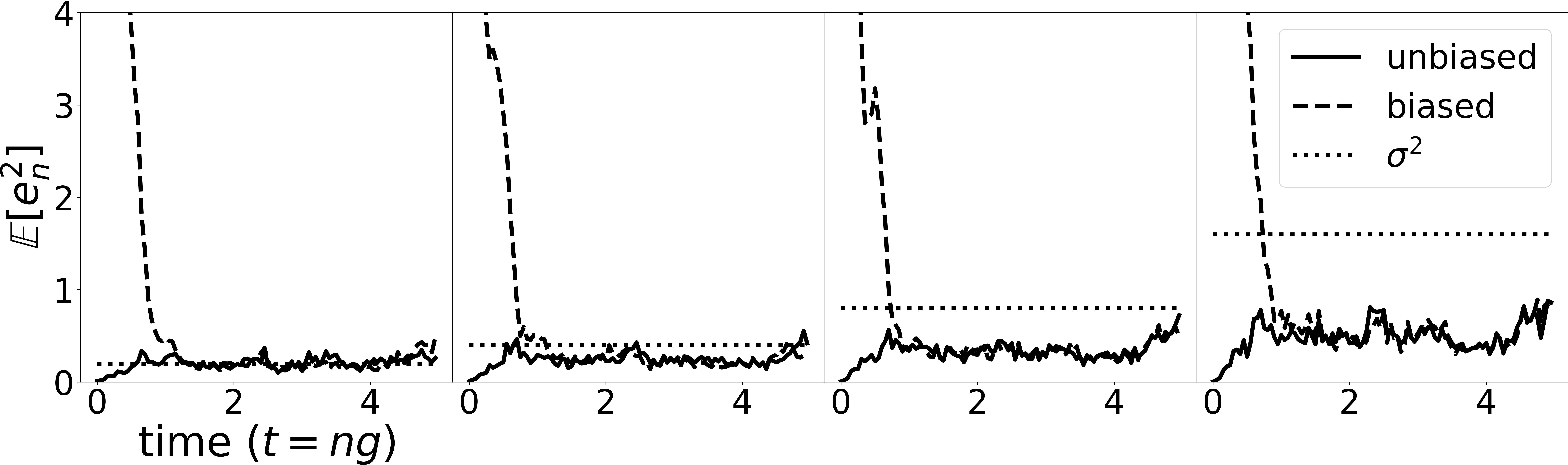} $\ \ $
    \includegraphics[width=0.48\textwidth]{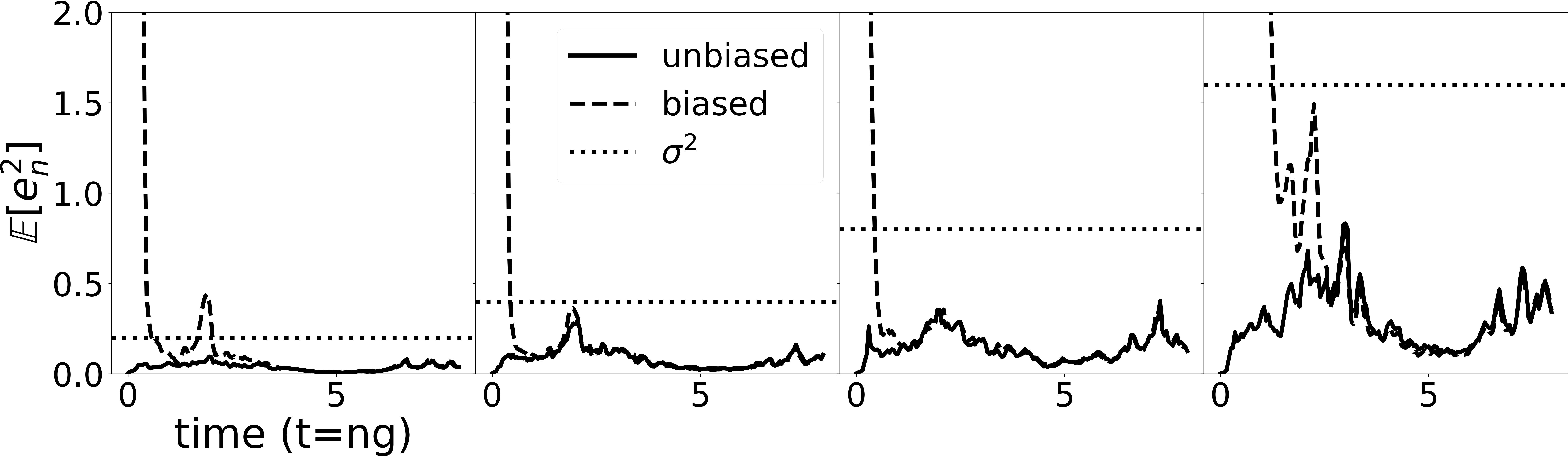}\\
   \includegraphics[width=0.48\columnwidth]{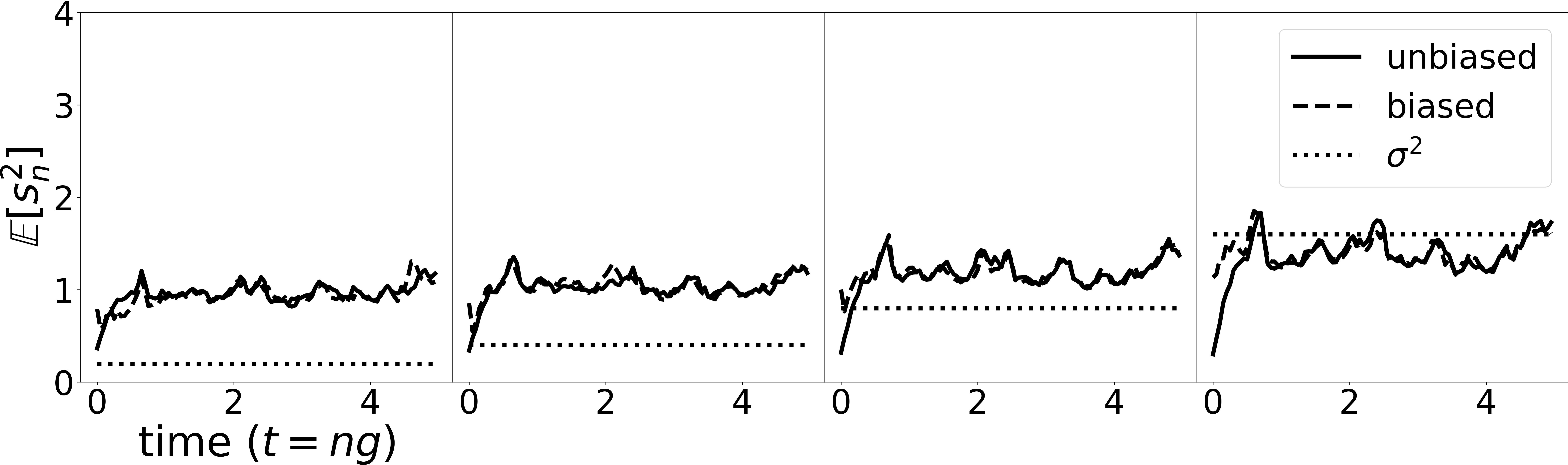} $\ \ $
    \includegraphics[width=0.48\textwidth]{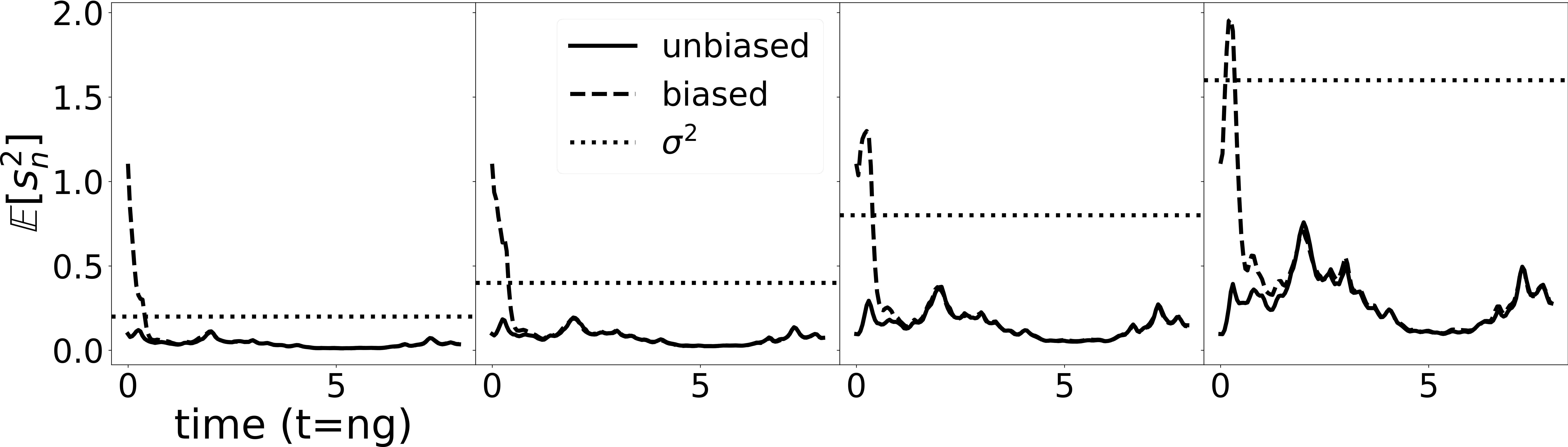}\\
   \includegraphics[width=0.48\columnwidth]{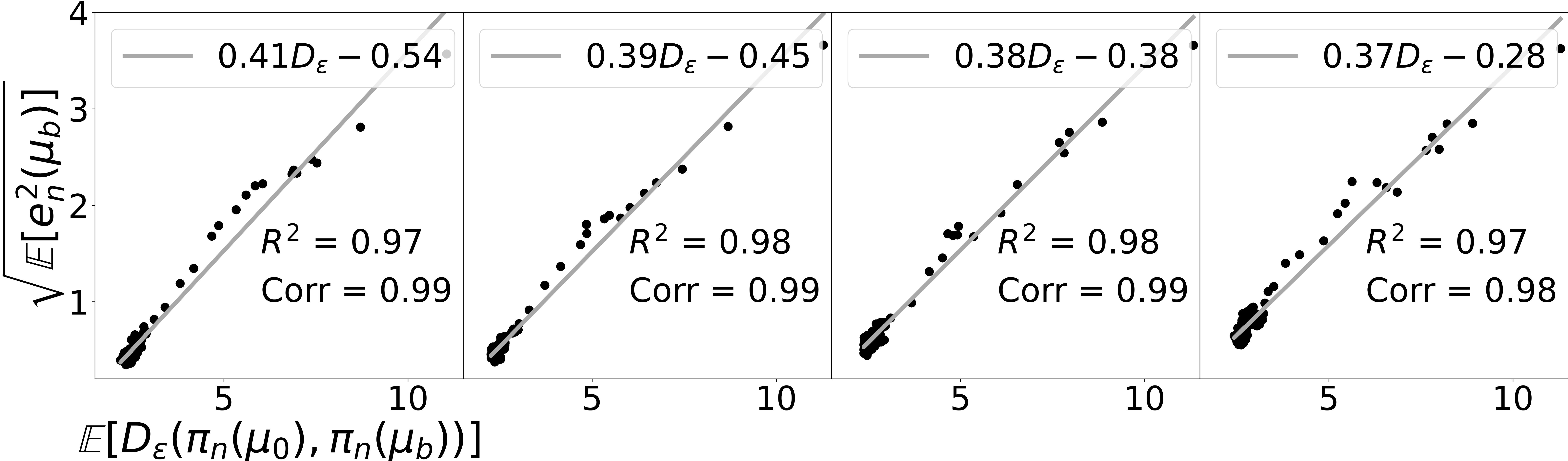} $\ \ $
    \includegraphics[width=0.48\textwidth]{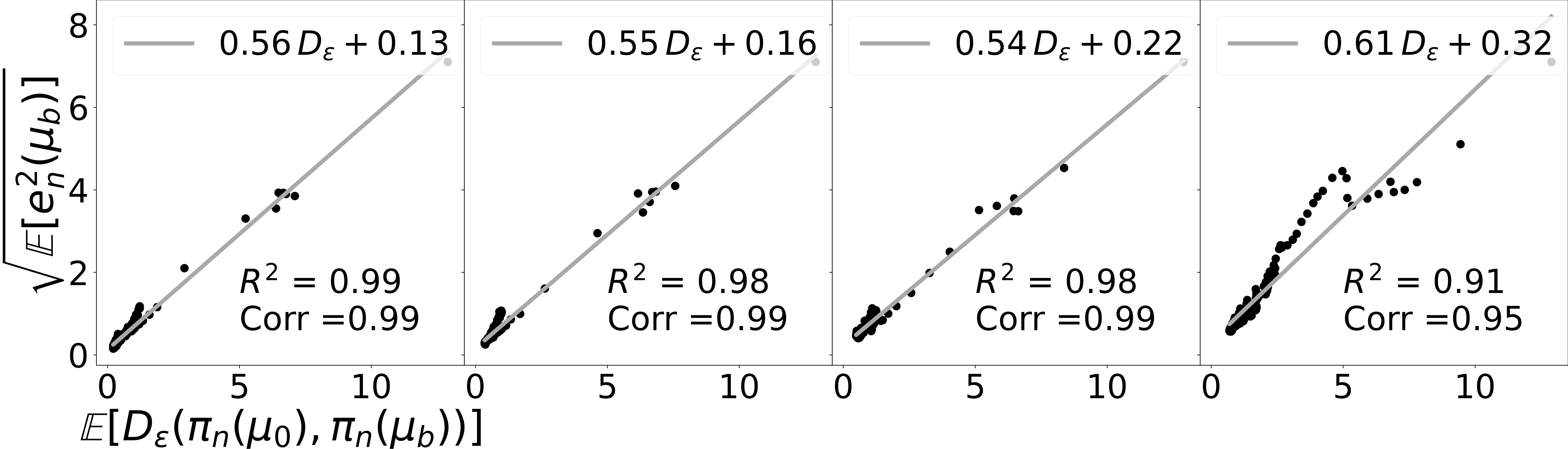}
\caption{Same as in figure~\ref{fig:bpf-enkf-fixed-ocov} with the left and right panels showing the results for PF and EnKF respectively, but with fixed time between observations $g = 0.05$, and each column containing the results for different observational error variances of $\sigma^2 = 0.2, \, 0.4, \, 0.8, \, 1.6$.}
\label{fig:bpf-enkf-fixed-ogap}
\end{figure}

We now discuss the results of numerical experiments with fixed observation gap $g = 0.05$, with varying observation covariance $\sigma^2 = 0.2, 0.4, 0.8, 1.6$. In figure \ref{fig:bpf-enkf-fixed-ogap}, we again note the exponential decrease of the distance $D_\varepsilon$ over time until it reaches a stationary value $c$. The parameter values obtained for the best-fit for different observation covariance are shown in table \ref{table:fixgap}.

\begin{table}[t!]
\centering
\begin{tabular}{|c|c|c|c|c|c|} 
 \hline
 
\multicolumn{2}{|c|}{$\sigma^2$} & $\bm{0.2}$ & $ \bm{0.4}$  & $\bm{0.8} $ & $\bm{1.6}$ \\ [0.5ex] 
\hline
\multirow{2}{*}{\text{a}} & \textbf{PF}& 7.842 $\pm$ 0.058 & 7.730 $\pm$ 0.046 & 8.153 $\pm$ 0.044 & 8.038 $\pm$ 0.048 \\\cline{2-6}
& \textbf{EnKF}& 10.61 $\pm$ 0.32 & 10.84 $\pm$ 0.30 & 10.69 $\pm$ 0.23 & 8.75 $\pm$ 0.22 \\
\hline
\multirow{2}{*}{$\lambda$}& \textbf{PF} & 2.442 $\pm$ 0.016 &  2.858 $\pm$ 0.017 & 2.859 $\pm$ 0.015 & 2.416 $\pm$ 0.012 \\ \cline{2-6}
& \textbf{EnKF} & 3.34 $\pm$ 0.16 &  3.70 $\pm$ 0.16 & 3.86 $\pm$ 0.14 & 1.507 $\pm$ 0.062 \\
\hline
\multirow{2}{*}{c} & \textbf{PF} & 2.4050 $\pm$ 0.0022 & 2.4144 $\pm$ 0.0015 & 2.5051 $\pm$ 0.0014 & 2.6554 $\pm$ 0.0019\\ \cline{2-6}
& \textbf{EnKF} & 0.470 $\pm$ 0.039 & 0.579 $\pm$ 0.035 & 0.838 $\pm$ 0.027 & 1.148 $\pm$ 0.041\\
\hline
\end{tabular}
\caption{Parameters of the best-fit for the mean $D_\varepsilon$ versus time as in~\eqref{eq:fit} with associated confidence intervals for fixed observation gap $g = 0.05$ and different observational error covariance $\sigma^2$ shown in the top row.}
\label{table:fixgap}
\end{table}

In contrast with the case of varying observational gap, the exponential rates for the PF stability are not affected by the change in observational uncertainty. While the rates for EnKF are again close to twice the Lyapunov exponent, the rates for PF are smaller.

The scaled $l_2$ error and the $D_\varepsilon$ achieve their stationary value around the same time as in the former case of fixed observation. As expected, this asymptotic value $c$ as well as the asymptotic values of the uncertainty $s_n$ and the bias $e_n$ all increase with increasing $\sigma^2$ for both PF and EnKF.  

We also see near perfect correlation in the scatter plots for the RMSE versus mean $D_\varepsilon$ as for both PF and EnKF, we get pearson correlation coefficient very close to 1. We remark that in our numerical experiments, either with varying observational time gap or with varying observational covariance, we did not notice any relation between stability and posterior uncertainty or precision, i.e., there did not seem to be any relation between $D_\varepsilon$ and $s_n$.

%% file: 06-appendix.tex
{\color{mypink}
\subsection{Definition~\ref{def-stab} implies definition \ref{def-stab-ra}}

We have introduced two different definitions of filter stability in section~\ref{ssec-def-stab}. Definition~\ref{def-stab} is used for the numerical computations in this paper while the other definition~\ref{def-stab-ra} is a commonly used definition using the notion of weak convergence. In this appendix, we prove that~\ref{def-stab} is stronger in the sense that it implies the other.
}
\begin{thm} Let $(M, d)$ be a compact metric space and $\{\alpha_n\}, \{\beta_n\}$ be sequences of random probability measures on $M$ with finite first and second moments and 
\begin{align}
    \lim_{n\to\infty}\mathbb E[D(\alpha_n, \beta_n)] = 0 \,,
\end{align}
where $D = W_1$ or $W_2$. Then for all globally Lipschitz functions $f: M \to \mathbb R$ we have
\begin{align}
    \lim_{n\to\infty}\mathbb E\left[\left|\int_M f\,d\alpha_n -\int_M f\,d\beta_n\right|\right] = 0 \,.
\label{eq-fconv} \end{align}

\begin{proof}
Let $S=\{f:M\to\mathbb R:\rm{Lip}(f)\le 1\}$ be the set of Lipschitz functions with Lipschitz constant not greater than 1.
If $f\in S$, by Kantorovich-Rubinstein duality we have,
\begin{align}
    \mathbb E\left[\left|\int_M f\,d\alpha_n -\int_M f\,d\beta_n\right|\right] \le 
    \mathbb E\left[\sup_{g\in S}\left|\int_M g\,d\alpha_n -\int_M g\,d\beta_n\right|\right]=\mathbb E[W_1(\alpha_n, \beta_n)] \,,\label{eq:con-KR}
\end{align}
which proves the assertion for the case $D=W_1$ for $f \in S$. For $D=W_2$, recall that for two measures $\mu, \nu$ on $M$ with finite first and second moments, 
\begin{align}
    W_p(\mu, \nu)^p = \inf_{\mathbb S}\mathbb E[d(X, Y)^p] \,, \label{eq:Wp-defn}
\end{align}
where $\mathbb S = \{\pi: \pi \text{ is the joint distribution of } X, Y \text{ with marginals  } \mu, \nu \text{ respectively}\}$. Therefore,
\begin{align}
    W_1(\mu, \nu)^2 &= \left(\inf_{\mathbb S}\mathbb E[d(X, Y)]\right)^2=\inf_{\mathbb S}\mathbb E[d(X, Y)]^2 \le \inf_{\mathbb S}\mathbb E[d(X, Y)^2], \quad\text{[by Jensen's inequality]}\\
    &\le W_2(\mu, \nu)^2,\quad\text{[by } \eqref{eq:Wp-defn}] \label{eq:Wp-ineq}
\end{align}
Combining \eqref{eq:con-KR} and \eqref{eq:Wp-ineq}, we immediately see that the assertion is true for the case $D = W_2$ for $f \in S$.

Now pick a function $\tilde{f}: M\to\mathbb R$ with $\text{Lip} (\tilde f)\le K$ for some $K>0$. Define $f=\frac{\tilde{f}}{K}$. Clearly, $f\in S$ and therefore,
\begin{align}
    \lim_{n\to\infty}\mathbb E\left[\left|\int_M \tilde f\,d\alpha_n -\int_M \tilde f\,d\beta_n\right|\right] =K\lim_{n\to\infty}\mathbb E\left[\left|\int_M f\,d\alpha_n -\int_M f\,d\beta_n\right|\right] =0 \,, 
\end{align}\label{thm:lip}
\end{proof}
\end{thm}

Theorem~\ref{thm:lip} immediately yields the following stronger statement.
\begin{cor}Under identical conditions, \eqref{eq-fconv} holds for any continuous $f:M\to\mathbb R$.

\begin{proof}

Note that real-valued, continuous functions on $M$ are uniform limits of locally Lipschitz functions on $M$. Since $M$ is compact and locally Lipschitz functions on compact metric spaces are globally Lipschitz, real-valued, continuous  functions on $M$ are uniform limits of globally Lipschitz functions on $M$. 

Let $\alpha_n - \beta_n=\gamma_n$ for brevity and let $f$ be a real-valued, continuous function on $M$.
\begin{align}
    \left|\int_M f\,d\alpha_n - \int_M f\, d\beta_n\right| = \left|\int_M f\,d\gamma_n\right| = \left|\int_M f^+\,d\gamma_n^+ + \int_M f^-\,d\gamma_n^- - \int_M f^-\,d\gamma_n^+ - \int_M f^+\,d\gamma_n^-\right| \,,
\end{align}
where $+, -$ denote the positive and negative parts of $f$ and $\gamma_n$. We now consider the case when the last sum of four terms is positive. The other case when it is negative is very similar.

Pick sequences of non-negative globally Lipschitz functions $f^+_m\uparrow f^+$ and $f^-_m\uparrow f^-$. Then by using dominated convergence theorem and Fatou's Lemma, we see that
\begin{align}
    \mathbb E\left[\int_M f^+\,d\gamma_n^+\right] = \mathbb E\left[\lim_{m\to\infty}\int_M f^+_m\,d\gamma_n^+\right] \le \mathbb E\left[\liminf_{m\to\infty}\int_M f^+_m\,d\gamma_n^+\right] \le \liminf_{m\to\infty}\mathbb E\left[\int_M f^+_m\,d\gamma_n^+\right] \,.
    \label{eq:ineq1}
\end{align}
Using very similar arguments for each of the four terms and adding the four inequalities, we see that 
\begin{align}
   \mathbb E\left[\left|\int_M f\,d\alpha_n - \int_M f\,d\beta_n\right|\right] \le \liminf_{m\to\infty}\mathbb E\left[\int_M f^+_m\,d\gamma_n^+\int_M f^-_m\,d\gamma_n^--\int_M f^-_m\,d\gamma_n^+-\int_M f^+_m\,d\gamma_n^-\right] \,. \label{eq:ineq-combo}
\end{align}
Define, $f_m = f_m^+ - f_m^-$ which is globally Lipschitz since sum of globally Lipschitz functions are globally Lipschitz. Since the RHS of \eqref{eq:ineq-combo} is non-negative, it must be equal to $\liminf_{m\to\infty}\mathbb E\left[\left|\int_M f_m\,d\gamma_n\right|\right]$. Therefore,
\begin{align}
    \mathbb E\left[\left|\int_M f\,d\alpha_n - \int_M f\,d\beta_n\right|\right]\le\liminf_{m\to\infty}\mathbb E\left[\left|\int_M f_m\,d\alpha_n - \int_M f_m\,d\beta_n\right|\right] \,.
\end{align}
Hence, given $\varepsilon > 0, \,\exists\, m_0$ such that $\forall \, m > m_0$, 
\begin{align}
\mathbb E\left[\left|\int_M f_m\,d\alpha_n - \int_M f_m\,d\beta_n\right|\right] > \mathbb E\left[\left|\int_M f\,d\alpha_n - \int_M f\,d\beta_n\right|\right] - \varepsilon \,.
\end{align}
Using theorem~\ref{thm:lip}, we see that the limit $n \to \infty$ of LHS is $0$ which gives
\begin{align}
    & 0 \ge \lim_{n\to\infty}\mathbb E\left[\left|\int_M f\,d\alpha_n - \int_M f\,d\beta_n\right|\right] - \varepsilon \,.
\end{align}
Since $\varepsilon$ can be chosen arbitrarily, we have
\begin{align}
    \lim_{n\to\infty}\mathbb E\left[\left|\int_M f\,d\alpha_n - \int_M f\,d\beta_n\right|\right] = 0  \,.
\end{align}
\end{proof}
\label{cor:cts}
\end{cor}

\begin{thm} Assume the filtering distributions $\hat\pi_n(\mu),\hat\pi_n(\nu)$ are supported on a compact subset of $\mathbb R^d$ and $D=W_2$ or $W_1$. If \eqref{eq-stablaw} holds then so does \eqref{def-stab-sumith}.
\begin{proof}
Direct consequence of corollary~\ref{cor:cts}.
\end{proof}
\label{thm:strength}
\end{thm}
{\color{mypink}\subsection{Convergence in $D_\varepsilon$ and $W_2$}\label{ssec-dw}
\begin{thm} Suppose 
$\exists$ a sequence $p_k: \lim_{k\to\infty}p_k=0$ and $  a_{kn}:=\mathbb E[D_{p_k}(\alpha_n, \beta_n)]$ is monotone decreasing in $k$ or, 
$a_{kn} \ge a_{k+1,  n}\;\forall\; k, n$ where $\{\alpha_n\}, \{\beta_n\}$ are sequences of random measures. If for some $j$
\begin{align}
    \limsup_{n\to\infty}\mathbb E[D_{p_j}(\alpha_n, \beta_n)] \le \delta,
\end{align}
then
\begin{align}\limsup_{n\to\infty}\mathbb E[W_2 (\alpha_n, \beta_n)]\le\delta
\end{align}

\begin{proof}
Recall that $\lim_{k\to\infty}D_{p_k} (\alpha_n, \beta_n) = W_2 (\alpha_n, \beta_n)$, see for example, theorem $1$ in \cite{genevay2018learning}. By Fatou's lemma,
\begin{align}
    \limsup_{n\to\infty}\mathbb E[W_2 (\alpha_n, \beta_n) \le \limsup_{n\to\infty}\liminf_{k\to\infty}\mathbb E[D_{p_k}(\alpha_n, \beta_n))] \le  \limsup_{n\to\infty}\mathbb E[D_{p_j}(\alpha_n, \beta_n))]\le\delta
\end{align}
\end{proof}
\label{thm:w2-limsup-bound}
\end{thm}
Theorem ~\ref{thm:w2-limsup-bound} is helpful in discussing the convergence of measures in Wasserstein metric by looking at convergence in $D_\varepsilon$ for a fixed $\varepsilon$. The assumption about existence of a monotone decreasing subsequence is consistent with numerical experiments as shown in figure~\ref{fig:bpf-eps-all}. Although Sinkhorn divergence is zero when we are comparing two identical distributions, while comparing two different samples from the same distribution we would expect a non-zero divergence which is why theorem~\ref{thm:w2-limsup-bound}  is stated in terms of a non-zero bound $\delta$. We note that a more detailed discussion about the behaviour of this non-zero divergence can be found in section~IV.A in~\cite{mandal2021stability}. With increasing sample size this non-zero divergence approaches zero. For a detailed look at the theoretical sample-complexity of Sinkhorn divergence see chapter~$3$ of~\cite{genevay2019entropy}.
\begin{figure}
    \centering
    \includegraphics[scale=0.15]{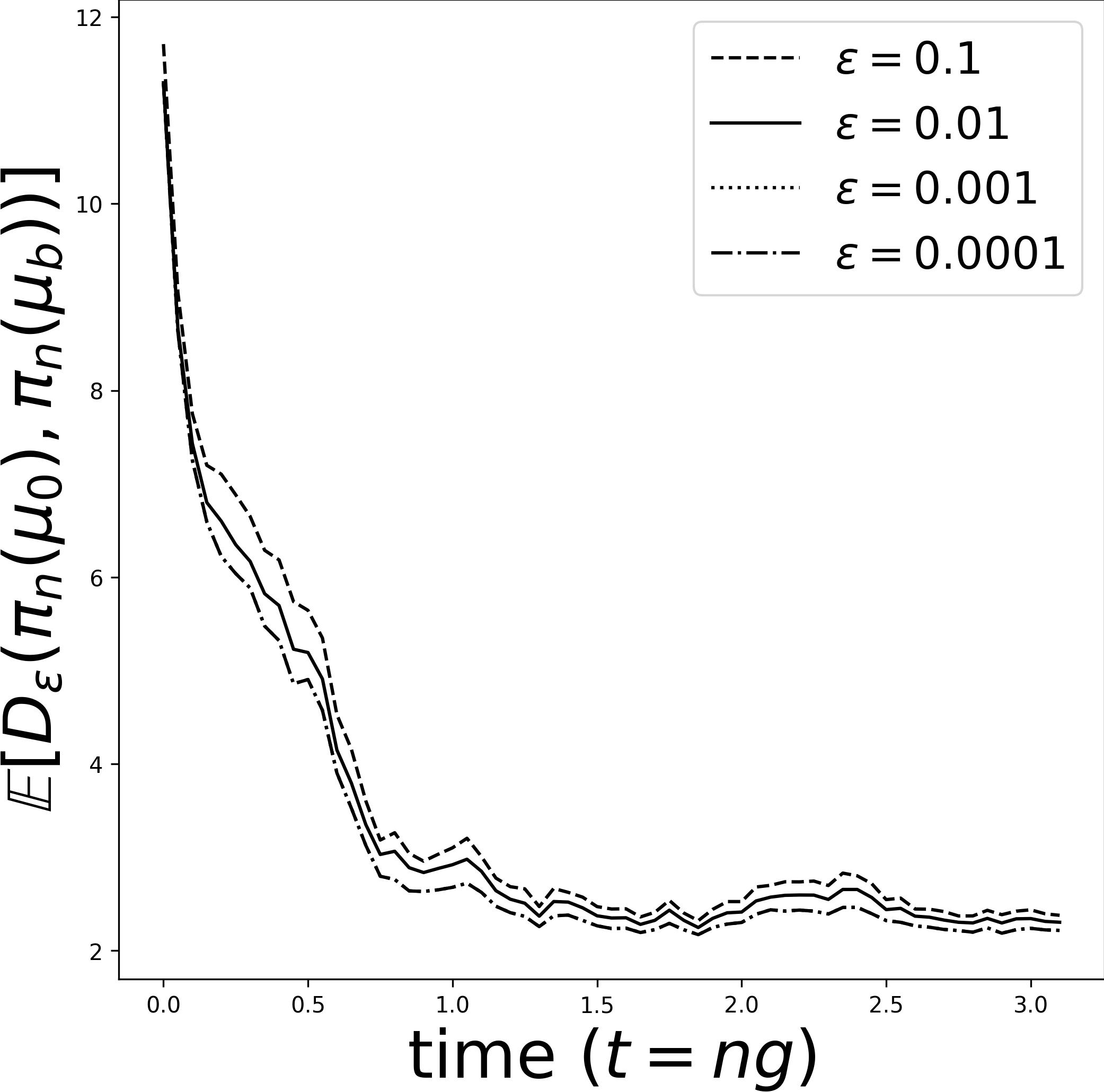}
    \caption{Change in mean $D_\varepsilon$ with $\varepsilon$. The smallest two  $\varepsilon$ values produce identical mean lines. The filtering distributions are generated by particle filter for observation gap  $=0.05$ and observation covariance $=0.4$.}
    \label{fig:bpf-eps-all}
\end{figure}
}

{\color{mypink}\subsection{Effect of varying the sample-size}\label{ssec-sample-size}
We find that averaging over $10$ observation realizations (for the expectation in~\eqref{eq-stablaw}) is sufficient for this study. To illustrate this, we repeated a representative experiment for a larger sample size  of $100$ and the resulting plots of the distance versus time for observation gap $g=0.05$ and $\sigma^2=0.4$ for the two filters are shown in figure~\ref{fig:obs-100}. Comparing with the plots in first row, second column of the particle filter and EnKF results in figure~\ref{fig:bpf-enkf-fixed-ogap}, we see that the results are qualitatively identical and quantitatively near-identical. Thus the choice of averaging over $10$ realizations suffices to capture the statistical features of the quantities we are studying.
\begin{figure}
    \includegraphics[scale=0.1]{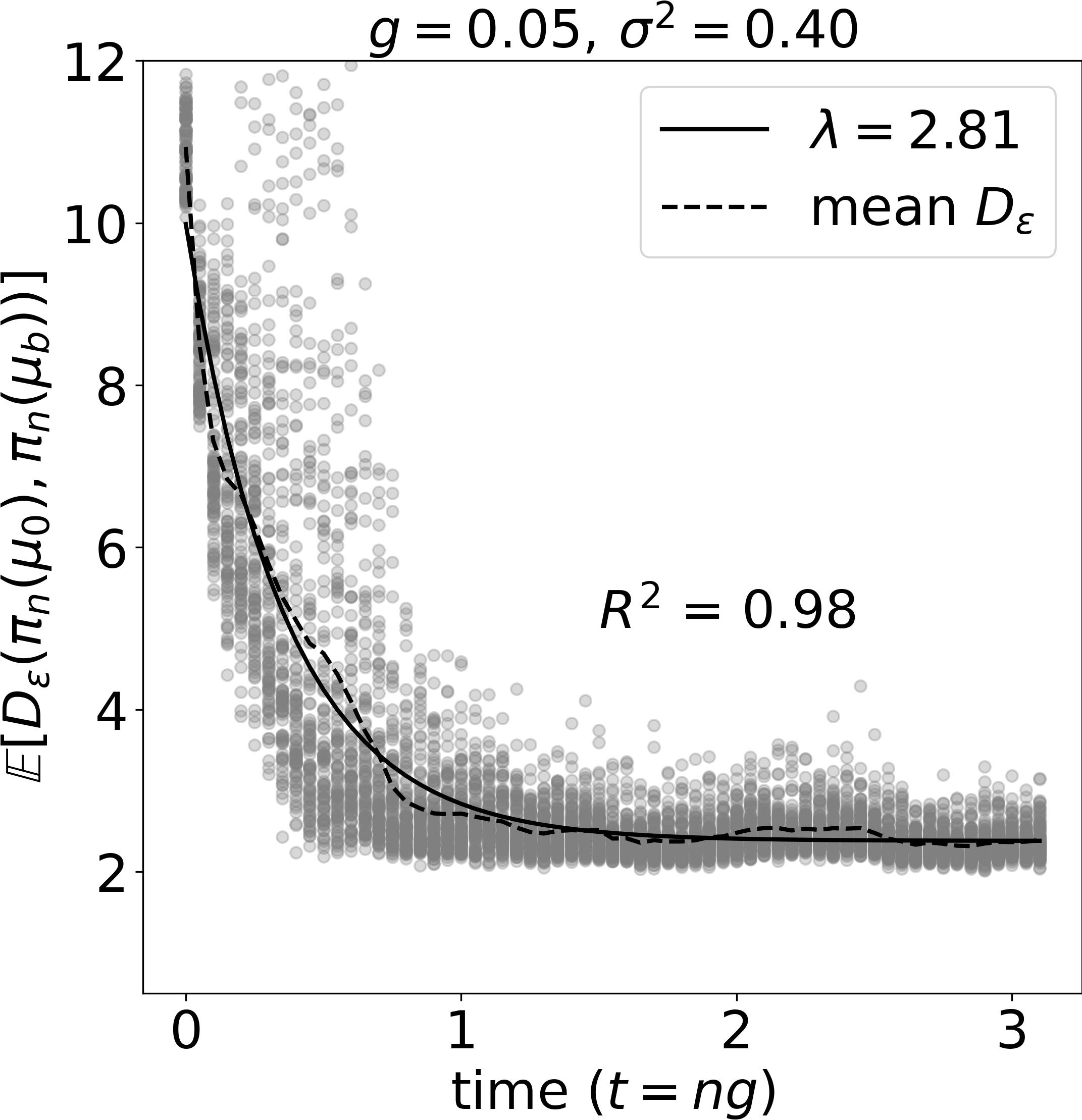}
    \includegraphics[scale=0.3]{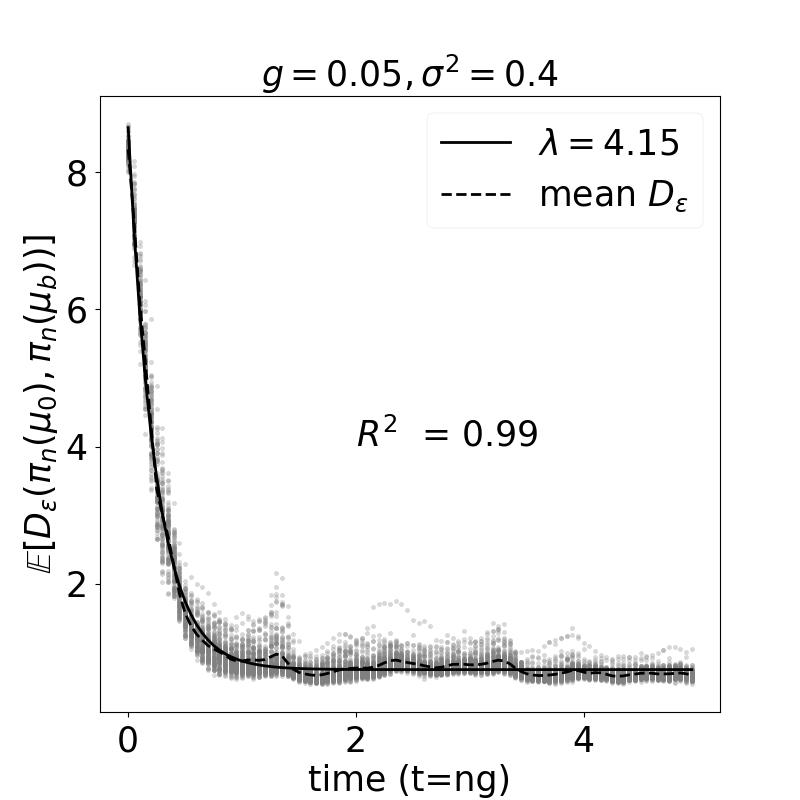}
    \caption{Mean $D_\varepsilon$ for $100$ observations for PF (left panel) and EnKF (right panel) with observation gap  $=0.05$ and observation covariance $=0.4$.}
    \label{fig:obs-100}
\end{figure}
}